\documentclass{amsart}
\usepackage[top=25mm, bottom=25mm, left=25mm, right=25mm]{geometry}
\usepackage{graphicx}
\usepackage{amsmath}
\allowdisplaybreaks[4]
\usepackage{amssymb}
\usepackage{amsthm}
\usepackage{amsfonts}
\usepackage{bm}
\usepackage{mathtools}
\usepackage{empheq}
\usepackage{tikz}

\usepackage{enumitem}
\usepackage{booktabs}

\usepackage{longtable}
\usepackage{float}

\theoremstyle{definition}

\newtheorem{definition}{Definition}[section]
\newtheorem{proposition}[definition]{Proposition}
\newtheorem{lemma}[definition]{Lemma}
\newtheorem{theorem}[definition]{Theorem}
\newtheorem{corollary}[definition]{Corollary}
\newtheorem{remark}[definition]{Remark}

\newtheorem{notation}[definition]{Notation}
\newtheorem{conjecture}[definition]{Conjecture}
\numberwithin{equation}{section}

\title{$p$-adic Congruences of Generalized Euler Numbers and Relations to Even Zeta Values}
\author{Yuta Nishibuchi}
\date{}

\address{Yuta Nishibuchi\\Mathematical Institute, Tohoku University\\6-3, Aoba, Aramaki, Aoba-ku, Sendai 980-8578, Japan}
\email{nishibuchi.yuta.s2@dc.tohoku.ac.jp}

\begin{document}

\begin{abstract}
    Generalized Euler numbers have previously been studied mainly from a combinatorial viewpoint. The purpose of this paper is to explore them from $p$-adic and analytic perspectives. To this end, we introduce congruential Euler numbers, a new family extending generalized Euler numbers. We establish several $p$-adic congruences for these numbers, including an answer to a conjecture related to Lehmer numbers. Furthermore, using complex analytic methods, we derive expressions for even zeta values in terms of congruential Euler numbers. These results suggest that certain types of these numbers may possess connections with arithmetic or analytic structures beyond their purely combinatorial role as generalizations of Euler-type numbers.

\end{abstract}

\maketitle

\tableofcontents

\section{Introduction}
The sequence of Euler numbers $\{E_n\}_{n=0}^\infty$ is defined by the exponential Taylor coefficients of the function
\[\text{sech }z=\frac{2}{e^z+e^{-z}}=\sum_{n=0}^\infty E_n \frac{z^n}{n!}.\]
In various areas of mathematics, several analogs and generalizations of Euler numbers have been introduced and studied.

In 1935, Lehmer defined the numbers $\{W_n\}$ by
\[\sum_{n=0}^{\infty}W_n\frac{z^n}{n!}=\frac{3}{e^z+e^{\omega z}+e^{\omega^2z}}=\left(\sum_{n=0}^{\infty}\frac{z^{3n}}{(3n)!}\right)^{-1}\]
where $\omega=-1/2+\sqrt{-3}/2$ is a cube root of unity. These numbers are called \textit{Lehmer numbers} today and are regarded as a combinatorially important analog of Euler numbers.

A common generalization of Euler numbers and Lehmer numbers is the \textit{generalized Euler numbers} $\{E_n^{(k)}\}$ defined by D. J. Leeming and R. A. MacLeod \cite{Leeming_and_MacLeod};
$$\sum_{n=0}^\infty E_n^{(k)}\frac{z^n}{n!}=\left(\sum_{n=0}^\infty \frac{z^{kn}}{(kn)!}\right)^{-1}$$
where $k\ge 2$ is a fixed integer. 
\footnote{They did not define $\{E_n^{(k)}\}$ in this way, but they proved that this is an equivalent definition in \cite[see Theorem 3.2]{Leeming_and_MacLeod}.} Regarding the combinatorics of generalized Euler numbers, B. E. Sagan provided a combinatorial interpretation of these numbers. 

An interesting property of these numbers is their various congruence equations. For example, in \cite{Gessel}, I. M. Gessel proved several congruences, including the conjectures of Leeming and MacLeod.
 In particular, T. Komatsu and G.-D. Liu investigated the $3$-adic congruence periodicities of Lehmer numbers in \cite{Komatsu_and_Liu} and proposed a conjecture;
\[3n\equiv 3m\pmod{2\cdot3^k}\rightarrow W_{3n}\equiv W_{3m}\pmod{3^{k+1}}\]
for any nonnegative integers $n,m$ and positive integer $k$.
This suggests $p$-adic properties of these numbers.

The purpose of this paper is to study $p$-adic congruences of generalized Euler numbers and demonstrate their usefulness from an analytic-number-theoretic perspective.

Motivated by the works of Barman and Komatsu \cite{Barman_and_Komatsu}, and T. Kameyama \cite{Kameyama}, in this paper, we define congruential Euler numbers $\{\mathcal{E}_{n}^{(N,j)}\}$ by
\[\sum_{n=0}^{\infty}\mathcal{E}_{n}^{(N,j)}\frac{z^n}{n!}=\left(\sum_{n=0}^{\infty}\frac{z^{Nn}}{(Nn+j)!}\right)^{-1}.\footnotemark\]
\footnotetext{The details of the background and motivation for this definition are discussed in Appendix A.} This is a further generalization of previously defined numbers. For example, the Euler numbers are the case of $(N,j)=(2,0)$, and the generalized Euler numbers are the case of $j=0$.

In this paper, we prove the $p$-adic congruences of congruential Euler numbers and establish relations between even zeta values and congruential Euler numbers of $(N,j)=(4,0),(4,2),(6,3)$. 

The main results are as follows. These results suggest that certain types of numbers, such as generalized Euler numbers, may have a relationship with other fields that goes beyond simply being a generalization of combinatorial objects.
\begin{theorem}
\label{MainThm}
    For any odd prime number $p$, positive integer $r$, integer $j$ with $0\le j\le p-1$, and nonnegative integer $n$,
    \begin{equation*}
    \mathcal{E}_{pn}^{(p,j)} + \mathcal{E}_{pn+p^r}^{(p,j)} \equiv 0 \mod{p^{r+\delta(j)}}\qquad
    \text{where }\delta(j)=
    \begin{cases}
        1 &(j=0)\\
        0 &(0<j<p).
    \end{cases}
\end{equation*}
\end{theorem}
This gives an affirmative answer to the conjecture of Komatsu and Liu \cite{Komatsu_and_Liu}.
\begin{proposition}
\label{Prop_for_conj_of_special_case}
    (1) For any positive integer $r$,
        \[\mathcal{E}_{4n+2^{r+1}}^{(4,0)}\equiv\mathcal{E}_{4n}^{(4,0)}\mod{2^r}\quad\text{ for any }n\ge 0.\]
        
    (2) For any positive integer $r$, there is a positive integer $n_0$ such that 
    \[\mathcal{E}_{6n+2\cdot3^r}^{(6,0)}\equiv\mathcal{E}_{6n}^{(6,0)}\mod{3^r}\quad\text{ for any }n\ge n_0.\]
\end{proposition}
\begin{theorem}
\label{MainThm2}
    For any positive integer $n$, the following equations hold.
    \begin{equation*}
        \lambda(4n)=\sum_{k=1}^\infty\frac{1}{(2k-1)^{4n}}=\frac{(-1)^{n+1}(\pi/\sqrt{2})^{4n}}{4(4n-1)!}\sum_{m=0}^{n-1}\binom{4n-1}{4m}\mathcal{E}_{4m}^{(4,0)},
    \end{equation*}
    \begin{equation*}
        \lambda(4n-2)=\sum_{k=1}^\infty\frac{1}{(2k-1)^{4n-2}}=\frac{(-1)^{n+1}(\pi/\sqrt{2})^{4n-2}}{4(4n-3)!}\sum_{m=0}^{n-1}\binom{4n-3}{4m}\mathcal{E}_{4m}^{(4,0)},
    \end{equation*}
    \begin{equation*}
        \zeta(4n)=\frac{(-1)^{n+1}(\sqrt{2}\pi)^{4n}}{4(4n-1)!(2^{4n}-1)}\sum_{m=0}^{n-1}\binom{4n-1}{4m}\mathcal{E}_{4m}^{(4,0)},
    \end{equation*}
    \begin{equation*}
        \zeta(4n-2)=\frac{(-1)^{n+1}(\sqrt{2}\pi)^{4n-2}}{4(4n-3)!(2^{4n-2}-1)}\sum_{m=0}^{n-1}\binom{4n-3}{4m}\mathcal{E}_{4m}^{(4,0)},
    \end{equation*}
    \begin{equation*}
        \zeta(4n)=\frac{(-1)^{n+1}(\sqrt{2}\pi)^{4n}}{4(4n+1)!}\sum_{m=0}^{n}\binom{4n+1}{4m}\mathcal{E}_{4m}^{(4,2)},
    \end{equation*}
    \begin{equation*}
        \zeta(4n-2)=\frac{(-1)^{n+1}(\sqrt{2}\pi)^{4n-2}}{4(4n-1)!}\sum_{m=0}^{n-1}\binom{4n-1}{4m}\mathcal{E}_{4m}^{(4,2)},
    \end{equation*}
    \begin{equation*}
        \zeta(6n)=\frac{(-1)^{n+1}(2\pi)^{6n}}{6(6n+2)!}\sum_{m=0}^{n}\binom{6n+2}{6m}\mathcal{E}_{6m}^{(6,3)},
    \end{equation*}
    \begin{equation*}
        \zeta(6n-4)=\frac{(-1)^{n+1}(2\pi)^{6n-4}}{6(6n-2)!}\sum_{m=0}^{n-1}\binom{6n-2}{6m}\mathcal{E}_{6m}^{(6,3)},
    \end{equation*}
    where $\zeta$ is the Riemann zeta function and
    $\lambda(s)=(1-2^{-s})\zeta(s)=L(s,\chi)$ (where $\chi$ is the trivial character modulo $2$) is the Dirichlet lambda function.
\end{theorem}

This paper is organized as follows. In Section 2, we introduce congruential Euler numbers and observe some of their basic properties. In Section 3, we discuss congruence relations of congruential Euler numbers; in Subsection 3.1, we prove Theorem \ref{MainThm}, and in Subsection 3.2, we propose a conjecture and give a proof of Proposition \ref{Prop_for_conj_of_special_case} as a special case of the conjecture. In Section 4, we prove Theorem \ref{MainThm2} using complex analysis. Appendix A provides an overview of related notions from 
prior work that serve as motivation for the definition of 
congruential Euler numbers. Appendix B offers computational evidence for the conjecture in Section 3.2.

\subsection{Acknowledgment}
I would like to thank Taiyo Kameyama for introducing the generalized Euler numbers and related topics, and for answering my questions. His preceding works helped me gain a quick understanding of the related field. I am profoundly grateful to my advisor, Nobuo Tsuzuki, for many valuable discussions and insightful advice, which greatly guided this research. I would also like to thank Kazuhiro Ito for participating in several discussions and for helpful comments. Lastly, I would like to thank Yasuo Ohno for comments on the manuscript.

\section{Congruential Euler numbers}
\begin{definition}
    For a fixed positive integer $N$ and a fixed nonnegative integer $j$ with $0\le j\le N-1$, we define the sequence of \textit{congruential Euler numbers} $\{\mathcal{E}_{n}^{(N,j)}\}_{n=0}^\infty$ by exponential Taylor coefficients of
    $$\sum_{n=0}^\infty\mathcal{E}_{n}^{(N,j)} \frac{z^n}{n!}=\left(\sum_{n=0}^\infty\frac{z^{Nn}}{(Nn+j)!}\right)^{-1}=\frac{Nz^j}{\sum_{k=0}^{N-1}\zeta_N^{-kj}\exp(\zeta_N^kz)}$$
    where $\zeta_N = e^{2\pi\sqrt{-1}/N}$ is a primitive $N$-th root of unity.
\end{definition}
\begin{remark}
\label{rem_generalization}
    Note that we have Euler numbers when we set $N=2$ and $j=0$ in the definition of congruential Euler numbers. Lehmer numbers are the case of $N=3$ and $j=0$. Setting $j=0$ with $N$ arbitrary gives generalized Euler numbers. We can naturally expand for the case $j\ge N$ by
    \[\sum_{n=0}^\infty\mathcal{E}_{n}^{(N,j)} \frac{z^n}{n!}=\left(\sum_{n=0}^\infty\frac{z^{Nn}}{(Nn+j)!}\right)^{-1}.\]
    Then, $N=1$ and $j=1$ provide well-known Bernoulli numbers $\{B_n\}_{n=0}^\infty$, which have the generating function
    $$\sum_{n=0}^{\infty} B_n \frac{z^n}{n!}=\frac{z}{e^z-1}.$$
\end{remark}

\begin{remark}
    By the definition of the congruential Euler numbers, $\mathcal{E}_{n}^{(N,j)}=0$ when $n$ is not a multiple of $N$. Thus, we mainly use the expressions $\sum\mathcal{E}_{Nn}^{(N,j)}z^{Nn}/(Nn)!$ in what follows.
\end{remark}

\begin{proposition}
\label{Prop_Fundamental_properties_j-Euler}
    For any integers $N,j$ with $N\ge 1$ and $j\ge 0$, 
    \begin{equation}
        \sum_{m=0}^{n}\binom{Nn+j}{Nm}\mathcal{E}_{Nm}^{(N,j)}=
        \begin{cases}
        j! & \text{if $n=0$,} \\
        0  & \text{if $n>0$,} 
        \end{cases}
    \end{equation}
\end{proposition}

\begin{proof}
    By definition, we have
    \begin{align*}
        &\sum_{n=0}^\infty\mathcal{E}_{Nn}^{(N,j)}\frac{z^{Nn}}{(Nn)!}\sum_{n=0}^\infty \frac{z^{Nn}}{(Nn+j)!}\\
        =&\sum_{n=0}^\infty\sum_{m=0}^{n}\frac{\mathcal{E}_{Nn}^{(N,j)}}{(Nm)!(N(n-m)+j)!}z^{Nn}\\ 
        =&\sum_{n=0}^\infty\left(\sum_{m=0}^{n}\binom{Nn+j}{Nm}\mathcal{E}_{Nm}^{(N,j)}\right)\frac{z^{Nn}}{(Nn+j)!}\\
        =&1.
    \end{align*}
    The conclusion is given by comparing the $n$-th coefficients.
\end{proof}

The next Corollary ensures that the numbers appearing in the congruences considered in this paper are p-adic integers. The first assertion is already known from \cite{Leeming_and_MacLeod}, but we include a proof here for later use.

\begin{corollary}
\label{Cor_denomintor_of_En}
    $\{\mathcal{E}_{Nn}^{(N,0)}\}_{n=0}^{\infty}$ are sequences of integers. Moreover, for any prime number $p$ and integer $j$ with $0\le j\le p-1$, $\{\mathcal{E}_{pn}^{(p,j)}\}_{n=0}^{\infty}$ are sequences of $p$-adic integers, i.e., the denominator (as an irreducible fraction) of $\mathcal{E}_{pn}^{(p,j)}$ does not have the prime factor $p$.
\end{corollary}
\begin{proof}
     The statement is given by the inductive use of Proposition \ref{Prop_Fundamental_properties_j-Euler}.
     We have 
     $$\mathcal{E}_{Nn}^{(N,j)}=-\left(\binom{Nn+j}{j}\right)^{-1}\sum_{m=0}^{n-1}\binom{Nn+j}{Nm}\mathcal{E}_{Nm}^{(N,j)}$$
     for $n\ge 1$.
     
     If $j=0$, then $\binom{Nn+j}{j}=1$, so $\mathcal{E}_{Nn}^{(N,j)}$ is an integer by induction.
     
     In the case that $N=p$ is a prime and $1\le j\le p-1$, $\binom{pn+j}{j}$ doesn't have a prime factor $p$, so $\mathcal{E}_{Nn}^{(N,j)}$ doesn't have a prime factor $p$ in its denominator by induction.
\end{proof}

\section{Congruence equations}
In this section, we will discuss several congruence equations regarding the congruential Euler numbers, and relate them to a conjecture from prior work and known results. Moreover, we propose a new conjecture and verify it in several special cases.

We start with a lemma.
\begin{lemma}
\label{Lem_binom_p_valuation}
    Let $p$ be a prime number and $v_p$ be the $p$-adic valuation.
    
    (1)Let $n,m$ be positive integers with $m\le n$. Then,
    \[v_p(\binom{n}{m})\ge v_p(n)-v_p(m).\]
    
    (2)Let $r$ be a positive integer and $m$ be a positive integer less than $p^r$. Then, 
    \[v_p(\binom{p^r}{m})=r-v_p(m).\]
\end{lemma}
\begin{proof}
    (1)Since $m>0$, we have
    \[\binom{n}{m}=\frac{n}{m}\binom{n-1}{m-1}.\]
    Taking the $p$-adic valuation on both sides and using the fact $v_p(\binom{n-1}{m-1})\ge0$, we obtain the conclusion.
    
    (2)Since $v_p(p^r-j)=v_p(j)$ for any integer $0<j<p^r$, we have
    \begin{align*}
        v_p(\binom{p^r}{m})&=v_p(\frac{p^r(p^r-1)\cdots(p^r-(n-1))}{m(m-1)\cdots2\cdot 1})\\
        &=\sum_{k=0}^{m-1}v_p(p^r-k)-\sum_{k=1}^{m}v_p(k)\\
        &=v_p(p^r)-v_p(m)\\
        &=r-v_p(m).
    \end{align*}
\end{proof}

For the convenience of the following discussions, we define the following functions.
\begin{notation}
\label{Not_F_H}
    For positive integer $N$ and $j=0,1,\cdots,N-1$, define functions $F_{N.j}(z)$ and $H_{N,j}(z)$ by
    $$F_{N,j}(z)=\sum_{n=0}^\infty \mathcal{E}_{Nn}^{(N,j)} \frac{z^{Nn}}{(Nn)!},\qquad  H_{N,j}(z)=\sum_{n=0}^\infty \frac{z^{Nn+j}}{(Nn+j)!}=\frac{1}{N}\sum_{k=0}^{N-1}\zeta_N^{-kj}\exp(\zeta_N^k z).$$
\end{notation}
\begin{remark}
    The following equations follow immediately from the definition.
    \begin{equation*}
        F_{N,j}(z)H_{N,j}(z)=z^j,\quad\quad
        \frac{d}{dz}H_{N,j}(z)=
        \begin{cases}
            H_{N,j-1}(z) & \text{if $j\ge 1$,} \\
            H_{N,N-1}(z) & \text{if $j=0$.} 
        \end{cases}
    \end{equation*}
    In particular, the $N$-th derivatives of $H_{N,j}(z)$ are themselves.
\end{remark}

\subsection{The case of odd primes}

First, we will prove the congruence equations $\mathcal{E}_{pn}^{(p,j)} + \mathcal{E}_{pn+p^r}^{(p,j)} \equiv 0 \mod{p^{r}}$ for odd primes $p$.

\begin{theorem}[Restatement of Theorem \ref{MainThm}]
    For any odd prime number $p$, positive integer $r$, integer $j$ with $0\le j\le p-1$, and nonnegative integer $n$,
    \begin{equation*}
    \mathcal{E}_{pn}^{(p,j)} + \mathcal{E}_{pn+p^r}^{(p,j)} \equiv 0 \mod{p^{r+\delta(j)}}\qquad
    \text{where }\delta(j)=
    \begin{cases}
        1 &(j=0)\\
        0 &(0<j<p).
    \end{cases}
\end{equation*}
\end{theorem}
By Corollary \ref{Cor_denomintor_of_En}, we can interpret this congruence in the $p$-adic integer ring $\mathbb{Z}_p$.
\begin{proof}
    In this proof, we abbreviated $F_{p,j}(z)$ as $F(z)$.
    Consider the $p^r$-th derivative of the definitional equation $F(z)H_{p,j}(z)=z^j$. Then, we have
    \begin{align*}
        &\sum_{k=0}^{p^r}\binom{p^r}{k}F^{(p^r-k)}(z)H_{p,j}^{(k)}(z)\\
        =&F(z)H_{p,j}(z)+\sum_{k_1=0}^{p^{r-1}-1}\sum_{k_0=0}^{p-1}\binom{p^r}{k_1p+k_0}F^{(p^r-k_1p-k_0)}(z)H_{p,j}^{(k_0)}(z)\\
        =&(F(z)+F^{(p^r)}(z))H_{p,j}(z)+\sum_{k_1=0}^{p^{r-1}-1}\sum_{k_0=1}^{p-1}\binom{p^r}{k_1p+k_0}F^{(p^r-k_1p-k_0)}(z)H_{p,j}^{(k_0)}(z)\\
        &+\sum_{k_1=1}^{p^{r-1}-1}\binom{p^r}{k_1p}F^{(p^r-k_1p)}(z)H_{p,j}(z)\\
        =&0.
    \end{align*}
    This gives
    $$\sum_{n=0}^\infty(\mathcal{E}_{pn}^{(p,j)} + \mathcal{E}_{pn+p^r}^{(p,j)})\frac{z^{pn}}{(pn)!}=-\frac{F(z)}{z^j}\sum_{k_1=0}^{p^{r-1}-1}\sum_{k_0=1}^{p-1}\binom{p^r}{k_1p+k_0}F^{(p^r-k_1p-k_0)}(z)H_{p,j}^{(k_0)}(z)-\sum_{k=1}^{p^{r-1}-1}\binom{p^r}{kp}F^{(p^r-kp)}(z)$$
    where we used the relationship $(H_{p,j}(z))^{-1}=F(z)/z^j$.

    Let $c(k_0)=1$ if $j<k_0$ and $c(k_0)=0$ if $j \ge k_0$.
    As for the double sum in the first term, we have
    \begin{align*}
        &\sum_{k_1=0}^{p^{r-1}-1}\sum_{k_0=1}^{p-1}\binom{p^r}{k_1p+k_0}F^{(p^r-k_1p-k_0)}(z)H_{p,j}^{(k_0)}(z)\\
        =&\sum_{k_1=0}^{p^{r-1}-1}\sum_{k_0=1}^{p-1}\binom{p^r}{k_1p+k_0}\left(\sum_{n=0}^\infty\mathcal{E}_{p(n-k_1)+p^r}^{(p,j)}\frac{z^{pn+k_0}}{(pn+k_0)!}\right)\left(\sum_{n=c(k_0)}^\infty\frac{z^{pn+j-k_0}}{(pn+j-k_0)!}\right)\\
        =&\sum_{k_1=0}^{p^{r-1}-1}\sum_{k_0=1}^{p-1}\binom{p^r}{k_1p+k_0}\left(\sum_{n=c(k_0)}^\infty\left(\sum_{m=0}^{n-c(k_0)}\binom{pn+j}{pm+k_0}\mathcal{E}_{p(m-k_1)+p^r}^{(p,j)}\right)\frac{z^{pn+j}}{(pn+j)!}\right)\\
        =&\sum_{n=1}^\infty\left(\sum_{k_1=0}^{p^{r-1}-1}\sum_{k_0=1}^{p-1}\binom{p^r}{k_1p+k_0}\sum_{m=0}^{n-c(k_0)}\binom{pn+j}{pm+k_0}\mathcal{E}_{p(m-k_1)+p^r}^{(p,j)}\right)\frac{z^{pn+j}}{(pn+j)!}\\
        &+\sum_{k_1=0}^{p^{r-1}-1}\sum_{k_0=1}^{j}\binom{p^r}{k_1p+k_0}\binom{j}{k_0}\mathcal{E}_{-pk_1+p^r}^{(p,j)}\frac{z^{j}}{j!}.
    \end{align*}
    By Lemma \ref{Lem_binom_p_valuation}, for any $0 \le k_1 \le p^{r-1}-1$ and $1 \le k_0 \le p-1$, $\binom{p^r}{k_1p+k_0}$ are multiples of $p^r$. Moreover, in the case of $j=0$, $\binom{pn+j}{pm+k_0}=\binom{pn}{pm+k_0}$ are multiples of $p$ for any $0 \le m \le n-1$ and $1\le k_0\le p-1$.
    Therefore, the first term is in the form
    \begin{align*}
        &\frac{F(z)}{z^j}\sum_{k_1=0}^{p^{r-1}-1}\sum_{k_0=1}^{p-1}\binom{p^r}{k_1p+k_0}F^{(p^r-k_1p-k_0)}(z)H_{k_0}(z)\\
        &=\frac{1}{z^j}\left(\sum_{n=0}^\infty\mathcal{E}_{pn}^{(p,j)}\frac{z^{pn}}{(pn)!}\right)\left(\sum_{n=0}^\infty(\text{multiple of }p^{r+\delta(j)})\frac{z^{pn+j}}{(pn+j)!}\right)\quad\left(\delta(j)=\begin{cases}
            1&(\text{if }j=0)\\
            0&(\text{if }j>0)
        \end{cases}\right)\\
        &=\frac{1}{z^j}\sum_{n=0}^\infty\left(\sum_{m=0}^n\binom{pn+j}{pm}\mathcal{E}_{pm}^{(p,j)}(\text{multiple of }p^{r+\delta(j)})\right)\frac{z^{pn+j}}{(pn+j)!}\\
        &=\frac{1}{z^j}\sum_{n=0}^\infty(\text{multiple of }p^{r+\delta(j)})\frac{z^{pn+j}}{(pn+j)!}\quad(\text{by Corollary \ref{Cor_denomintor_of_En}})\\
        &=\sum_{n=0}^\infty\frac{(\text{multiple of }p^{r+\delta(j)})}{(pn+j)\cdots (pn+1)}\frac{z^{pn}}{(pn)!}\\
        &=\sum_{n=0}^\infty(\text{multiple of }p^{r+\delta(j)})\frac{z^{pn}}{(pn)!}
    \end{align*}
    where ``multiple of $p^{r+\delta(j)}$" means some rational number whose numerator is a multiple of $p^{r+\delta(j)}$ and a denominator is not a multiple of $p$.
    
    As for the second term,
    \begin{align*}
        \sum_{k=1}^{p^{r-1}-1}\binom{p^r}{kp}F^{(p^r-kp)}(z)&=\sum_{k=1}^{p^{r-1}-1}\binom{p^r}{kp}F^{(kp)}(z)\quad(\text{replaced }k\text{ by }p^{r-1}-k)\\
        &=\sum_{k=1}^{p^{r-1}-1}\binom{p^r}{kp}\left(\sum_{n=0}^{\infty}\mathcal{E}_{p(n+k)}^{(p,j)}\frac{z^{pn}}{(pn)!}\right)\\
        &=\sum_{n=0}^{\infty}\left(\sum_{k=1}^{p^{r-1}-1}\binom{p^r}{kp}\mathcal{E}_{p(n+k)}^{(p,j)}\right)\frac{z^{pn}}{(pn)!}\\
        &=\sum_{n=0}^{\infty}\left(\sum_{k=1}^{(p^{r-1}-1)/2}\binom{p^r}{kp}(\mathcal{E}_{p(n+k)}^{(p,j)}+\mathcal{E}_{p(n+p^{r-1}-k)}^{(p,j)})\right)\frac{z^{pn}}{(pn)!}\quad(p\text{ is odd}).
    \end{align*}
    If $r=1$, this is an empty sum and equals 0.
    In the case of $r\ge 2$, we can show that $p^{r+\delta(j)}$ divides $\binom{p^r}{kp}(\mathcal{E}_{p(n+k)}^{(p,j)}+\mathcal{E}_{p(n+p^{r-1}-k)}^{(p,j)})$ for all $1 \le k \le p^{r-1}-1$ by induction. 
    
    Indeed, for each fixed $k$, we can write $k=a_ip^{i}+\cdots +a_{r-2}p^{r-2}$ with $0\le i\le r-2$, $0\le a_i,\cdots, a_{r-2}\le p-1$, and $a_i\ne 0$. Then $\binom{p^r}{kp}$ is a multiple of $p^{r-i-1}$ by Lemma \ref{Lem_binom_p_valuation} and $$\mathcal{E}_{p(n+k)}^{(p,j)}+\mathcal{E}_{p(n+p^{r-1}-k)}^{(p,j)}\equiv (-1)^{M}(\mathcal{E}_{p(n)}^{(p,j)}+\mathcal{E}_{p(n+p^{r-1})}^{(p,j)})\equiv 0\mod{p^{i+1+\delta(j)}}(\exists M)$$ by the inductive hypothesis.

    In conclusion, we have
    $$\sum_{n=0}^\infty(\mathcal{E}_{pn}^{(p,j)} + \mathcal{E}_{pn+p^r}^{(p,j)})\frac{z^{pn}}{(pn)!}=\sum_{n=0}^\infty(\text{multiple of }p^{r+\delta(j)})\frac{z^{pn}}{(pn)!}.$$
\end{proof}

This theorem provides answers to the conjecture proposed by Komatsu and Liu \cite[Conjecture 1]{Komatsu_and_Liu}.
To discuss its application to the conjecture, we recall the definition of Lehmer numbers.

\begin{definition}
    The Lehmer numbers $\{W_n\}_{n=0}^{\infty}$ are the exponential Taylor coefficients of the generating function
    \[\frac{3}{e^z+e^{\omega z}+e^{\omega^2z}}=\sum_{n=0}^{\infty}W_n\frac{z^n}{n!}\]
    where $\omega=-1/2+\sqrt{-3}/2$ is a cube root of unity.
\end{definition}
As referred to in Remark \ref{rem_generalization}, this is the special case of congruential Euler numbers when we set $N=3$ and $j=0$.

\begin{corollary}[Conjecture 1 in \cite{Komatsu_and_Liu}]
    
    Let $\{W_n\}$ be Lehmer numbers.
    For any nonnegative integers $n,m$ and positive integer $k$, if 
    \[3n\equiv 3m\pmod{2\cdot 3^k},\]
    then
    \[W_{3n}^{(3,0)}\equiv W_{3m}^{(3,0)}\pmod{3^{k+1}}\]
\end{corollary}
\begin{proof}
    Since the Lehmer numbers are the congruential numbers of $(N.j)=(3,0)$, by Theorem \ref{MainThm},
    \[W_{3n}\equiv-W_{3n+3^k}\equiv W_{3n+2\cdot 3^{k}}\pmod{3^{k+1}}\]
    for any nonnegative integer $n$. This is equivalent to the statement.
\end{proof}

The Theorem \ref{MainThm} gives the congruence properties of the case where $N=p$ is an odd prime. In the latter half of this section, we discuss congruence properties for the case where $N$ is a composite number.

If $N$ is a power of an odd prime and $j=0$, then, combined with Gessel's result, we have a similar congruence for small $r$.
In \cite{Gessel}, I. M. Gessel showed useful congruence properties for the generalized Euler numbers, which is the congruential Euler numbers in the case $j=0$. Using our notations, the statement is as follows.
\begin{proposition}[Gessel, Theorem 2.1 in \cite{Gessel}]

    For any prime $p$, any nonnegative integer $n$, and positive integer $m,r$,
    \[\mathcal{E}_{p^kmn}^{(p^km,0)}\equiv\mathcal{E}_{p^{k-1}mn}^{(p^{k-1}m,0)}\mod{p^{3k-\varepsilon}}\]
    where $\varepsilon=1$ for $p=2,3$ and $\varepsilon=0$ for $p\ge 5$.
\end{proposition}

By using Gessel's result, we have the following:
\begin{proposition}
\label{Prop_p_power}
    Let $p$ be an odd prime number, and let $\varepsilon=1$ if $p=3$, and $\varepsilon=0$ if $p\ge5$.
    Then, for any integer $1\le r\le 5-\varepsilon$, positive integer $k$, and nonnegative integer $n$,
    \[\mathcal{E}_{p^k(n+p^{r-1})}^{(p^k,0)}+\mathcal{E}_{p^kn}^{(p^k,0)}\equiv0 \mod{p^{r+1}}.\]
\end{proposition}
\begin{proof}
    Gessel's result states that $\cdots\equiv\mathcal{E}_{p^kn}^{(p^k,0)}\equiv\cdots\equiv\mathcal{E}_{p^2n}^{(p^2,0)}\equiv\mathcal{E}_{pn}^{(p,0)}\pmod{p^{6-\varepsilon}}$. So the congruence of Theorem \ref{MainThm} propagates to the case of $N=p^k$ if $r\le 5-\varepsilon$.
\end{proof}

\subsection{The case of composite numbers}

By Proposition \ref{Prop_p_power}, there seems $p$-adic congruences for some composite number $N$. In fact, computer calculations imply similar periodicities of remainders for $N$ of the form $(\text{a divisor of }p-1)\times p$. 

\begin{conjecture}
\label{conjecture}
    Let $p$ be an prime number and $m$ be a divisor of $p-1$ or $m=2$. Let $q$ be a least common multiple of $2$ and $p-1$. 
    
    Then, $\mathcal{E}_{mpn}^{(mp,j)}$ are $p$-adic integers and for any positive integer $r$, there is a positive integer $n_0$ such that 
    \[\mathcal{E}_{mpn+qp^r}^{(mp,j)}\equiv\mathcal{E}_{mpn}^{(mp,j)}\mod{p^r}\quad\text{ for any }n\ge n_0.\]
\end{conjecture}

For the results of numerical calculations that support this conjecture, see Appendix B. The calculations suggest $n_0$ can be taken as less than or equal to $r$.

Here, we provide a proof for a special case of the conjecture \ref{conjecture} as further evidence of it.

\begin{proposition}[Restatement of Proposition \ref{Prop_for_conj_of_special_case}]
    (1) For any positive integer $r$,
        \[\mathcal{E}_{4n+2^{r+1}}^{(4,0)}\equiv\mathcal{E}_{4n}^{(4,0)}\mod{2^r}\quad\text{ for any }n\ge 0.\]
    (This is the special case of the Conjecture \ref{conjecture}; $p=2,m=2,j=0$.)
        
    (2) For any positive integer $r$, there is a positive integer $n_0$ such that 
    \[\mathcal{E}_{6n+2\cdot3^r}^{(6,0)}\equiv\mathcal{E}_{6n}^{(6,0)}\mod{3^r}\quad\text{ for any }n\ge n_0.\]
    (This is the special case of the Conjecture \ref{conjecture}; $p=3,m=2,j=0$.)
\end{proposition}
\begin{notation}
     In the following proof of Proposition \ref{Prop_for_conj_of_special_case} and its auxiliary lemmas \ref{Lem_for_prf_for_conj_of_special_case} and \ref{Lem_2_for_special_case}, for any formal series $f(z)=\sum_{n=0}^{\infty}a_nz^n/n!$ and $g(z)=\sum_{n=0}^{\infty}b_nz^n/n!$ with $a_n,b_n\in\mathbb{Z}_p$,
    \[\text{`` }f\equiv g\mod{3^r}\text{ "}\quad\text{means}\quad a_n\equiv b_n\mod{3^r}\text{ for all }n\ge 0.\]
    Note that $H_{N,0}$ and, by Corollary \ref{Cor_denomintor_of_En}, $F_{N,0}$ have integer coefficients.
\end{notation}
\begin{proof}[Proof of Proposition \ref{Prop_for_conj_of_special_case}]
    Let $p=2$ or $3$. In this proof, we will use the following abbreviations;
    \[F(z)=F_{2p,0}(z)=\sum_{n=0}^{\infty}\mathcal{E}_{2pn}^{(2p,0)}\frac{z^{2pn}}{(2pn)!},\quad H(z)=H_{2p,0}(z)=\sum_{n=0}^{\infty}\frac{z^{2pn}}{(2pn)!}.\]

    For each integer $N\ge0$ let's define formal series $D_N$ to satisfy
    \[\frac{d^N}{dz^N}\left(\frac{1}{H}\right)=(-1)^N\frac{D_N}{H}.\]
    We can write $D_N$ explicitly by Faà di Bruno's formula; $D_N=\det(\binom{i}{j-1}H^{(i-j+1)})_{1\le i,j\le N}/H^N$ ($\binom{i}{j-1}=0$ if $i<j-1$). However, in this discussion, we do not need this explicit expressions. What we use here is the fact that $D_N$ is also of the form $\sum_{n=0}^{\infty} a_nz^n/n!,$ where $a_n\in\mathbb{Z}_p$.

    Then, by the definition of $F,D_N$ and by Lemma \ref{Lem_for_prf_for_conj_of_special_case} below, 
    \begin{align*}
        F-F^{(2\cdot p^r)}&=F-\frac{d^{2\cdot p^r}}{dz^{2\cdot p^r}}\left(\frac{1}{H}\right)\\
        &=F-\frac{D_{2\cdot p^r}}{H}\\
        &=F\left(1-D_{2\cdot p^r}\right)\\
        &\equiv 
        \begin{cases}
            F\cdot0\quad(\text{if }p=2)\\
            F\frac{(H^{(3)}/H)^2-1}{1+3(H^{(3)}/H)^2}\quad(\text{if }p=3)\\
        \end{cases}
    \mod{p^r}.\quad(\text{by Lemma \ref{Lem_for_prf_for_conj_of_special_case}})
    \end{align*}
    So this finishes the proof for the case of $p=2$. 
    
    As for the case of $p=3$, we have
    \begin{align*}
         F-F^{(2\cdot 3^r)}&=F\frac{(H^{(3)})^2-H^2}{H^2+3(H^{(3)})^2}\\
        &=((H^{(3)})^2-H^2)\frac{1}{H^3+3H(H^{(3)})^2}\mod{3^r}.
    \end{align*}

    In Lemma \ref{Lem_2_for_special_case}, we will show that both $(H^{(3)})^2-H^2\text{ and }1/(H^3+3H(H^{(3)})^2)$ are congruent modulo $3^r$ to some polynomial depending on $r$.
    So we obtain the following;
    \[F-F^{(2\cdot3^r)}\equiv(\text{polynomial depending on $r$})\mod{3^r}.\]
    This means that for each $r$, there exists some $n_0$ such that $\mathcal{E}_{6n+2\cdot3^r}^{(6,0)}\equiv\mathcal{E}_{6n}^{(6,0)}\mod{3^r}$ for all $n\ge n_0$.
\end{proof}

\begin{lemma}
\label{Lem_for_prf_for_conj_of_special_case}
    Let $p=2$ or $3$.
    Let $H(z)=\sum_{n=0}^{\infty}z^{2pn}/(2pn)!$  and 
    $D_N$ be as defined above. 

    Let $T(z)\coloneqq H^{(p)}(z)/H(z)$ and $v_p$ be the $p$-adic valuations.

    (1) If $p=2$, for any $m>0$,
    \begin{equation}
    \label{Equ_Xm_2}
        X_m\coloneqq 1-D_{2m}\equiv
        \begin{cases}
            2^{v_2(m)}(1+T^2) \quad(\text{if }m\text{ : even})\\
            1+T \quad(\text{if }m\text{ : odd})\\
        \end{cases}
        \quad\pmod{2^{v_2(2m)}}.
    \end{equation}
    
    (2) If $p=3$, for any $m>0$,
    \begin{equation}
    \label{Equ_Xm_3}
        X_m\coloneqq (-1)^m\left(1-D_{3m}\right)\equiv
        \begin{cases}
            \displaystyle\frac{(1+T)(T-1)}{1+3T^2} \quad(\text{if }m\text{ : even})\\
            \displaystyle\frac{(1-3T)(T-1)}{1+3T^2} \quad(\text{if }m\text{ : odd})\\
        \end{cases}
        \quad\pmod{3^{v_3(3m)}}.
    \end{equation}
    Note that \[\frac{1}{1+3T^2}=1-3T^2+9T^4-\cdots\] since the constant term of $T$ is $0$.
\end{lemma}
\begin{proof}
    We also define $X_0=0$, consistently with $D_0=1$.
    Taking the $N$-th derivative of $H\cdot(1/H)=1$, we have
    \[\sum_{k=0}^{N}\binom{N}{k}H^{(N-k)}\frac{d^k}{dz^k}\left(\frac{1}{H}\right)=\sum_{k=0}^{N}\binom{N}{k}(-1)^{k}\frac{H^{(N-k)}}{H}D_k=0\quad(N>0).\]
    Thus,
    \begin{equation}
    \label{Equ_not _congru}
        \sum_{k=0}^{N}\binom{N}{k}(-1)^{k}\frac{H^{(N-k)}}{H}\left(1-D_k\right)=\sum_{k=0}^{N}\binom{N}{k}(-1)^k\frac{H^{(N-k)}}{H}\quad(N>0).
    \end{equation}

    \underline{In the case of $p=2$}
    
    Let substitute $N=2m\text{ }(m>0)$ for (\ref{Equ_not _congru}). If $k$ is odd, $\binom{2m}{k}\equiv 0\mod{2^{v_2(2m)}}$ by Lemma \ref{Lem_binom_p_valuation}. So the left side hand of (\ref{Equ_not _congru}) will be
    \[\sum_{k=0}^{2m}\binom{2m}{k}(-1)^{k}\frac{H^{(2m-k)}}{H}\left(1-D_k\right)\equiv\sum_{k=0}^{m}\binom{2m}{2k}\frac{H^{(2m-2k)}}{H}X_m\mod{2^{v_2(2m)}}.\]
    Similarly, right side hand of (\ref{Equ_not _congru}) will be
    \begin{align*}
        \sum_{k=0}^{2m}\binom{2m}{k}(-1)^{k}\frac{H^{(2m-k)}}{H}&\equiv\sum_{k=0}^{m}\binom{2m}{2k}\frac{H^{(2m-2k)}}{H}\mod{2^{v_2(2m)}}\\
        &=
        \begin{cases}
            \displaystyle\sum_{\substack{k=0\\k\text{:even}}}^{m}\binom{2m}{2k}+\frac{H^{(2)}}{H}\sum_{\substack{k=0\\k\text{:odd}}}^{m}\binom{2m}{2k}\quad(\text{if $m>0$ : even})\\
            \displaystyle\frac{H^{(2)}}{H}\sum_{\substack{k=0\\k\text{:even}}}^{m}\binom{2m}{2k}+\sum_{\substack{k=0\\k\text{:odd}}}^{m}\binom{2m}{2k}\quad(\text{if $m>0$ : odd})\\
        \end{cases}\mod{2^{v_2(2m)}}.\\
    \end{align*}
    The sums of binomial coefficients can be computed explicitly. First,
    \begin{align*}
        \sum_{\substack{k=0\\k\text{:even}}}^m\binom{2m}{2k}&=\sum_{\substack{0\le k\le 2m\\k\equiv0\,(\mathrm{mod}\,4)}}\binom{2m}{k}\\
        &=\sum_{k=0}^{2m}\binom{2m}{k}\frac{1}{4}\sum_{j=0}^{3}(\sqrt{-1})^{jk}\\
        &=\frac{1}{4}\sum_{j=0}^{3}(1+(\sqrt{-1})^j)^{2m}\\
        &=\frac{1}{4}(2^{2m}+(2\sqrt{-1})^m+0^{2m}+(-2\sqrt{-1})^m)\\
        &=\begin{cases}
            2^{2m-2}+(-1)^{m/2}2^{m-1}&(\text{if $m$ : even})\\
            2^{2m-2}&(\text{if $m$ : odd})\\
        \end{cases}
    \end{align*}
    And second, since
    \[\sum_{\substack{k=0\\k\text{:even}}}^m\binom{2m}{2k}+\sum_{\substack{k=0\\k\text{:odd}}}^m\binom{2m}{2k}=\sum_{\substack{k=0\\k\text{:even}}}^m\binom{2m}{k}=\frac{1}{2}\sum_{\substack{k=0}}^m\binom{2m}{k}=2^{2m-1},\]
    we have
    \begin{equation*}
        \sum_{\substack{k=0\\k\text{:odd}}}^m\binom{2m}{2k}=
        \begin{cases}
            2^{2m-2}-(-1)^{m/2}2^{m-1}&(\text{if $m$ : even})\\
            2^{2m-2}&(\text{if $m$ : odd})\\
        \end{cases}
    \end{equation*}
    Now, for any $m\ge 3$, $v_2(2m)=v_2(m)+1\le m-1$. So the recurrence formula will be
    \begin{equation*}
        \sum_{k=0}^{m}\binom{2m}{2k}\frac{H^{(2m-2k)}}{H}X_m\equiv
        \begin{cases}
            1+T & (\text{if }m=1)\\
            2(1+T) & (\text{if }m=2)\\
            0 & (\text{if }m\ge 3)\\
        \end{cases}
        \mod{2^{v_2(2m)}}.
    \end{equation*}
    We can see the equation (\ref{Equ_Xm_2}) holds for $m=1$ by this, and thus $\binom{4}{2}X_1\equiv\binom{4}{2}(1+T)\mod{2^2}$ since $v_2(\binom{4}{2})=1$. Therefore,
    \[X_2\equiv2(1+T)-\binom{4}{2}TX_1\equiv2(1+T^2)\mod{2^2}.\]
    So (\ref{Equ_Xm_2}) also holds for $m=2$.
    Assume that (\ref{Equ_Xm_2}) holds for all $1\le k <m\text{ }(m\ge3)$. Then, since $v_2(\binom{2m}{2k})\ge v_2(2m)-v_2(2k)$ by Lemma \ref{Lem_binom_p_valuation}, 
    \begin{equation}
    \label{Eq_asuump_2}
        \binom{2m}{2k}X_k\equiv
        \begin{cases}
            \displaystyle\binom{2m}{2k}2^{v_2(k)}(1+T^2) &(\text{if }k\ge1\text{ : even})\\
            \displaystyle\binom{2m}{2k}(1+T) &(\text{if }k\ge1\text{ : odd})\\
        \end{cases}
        \quad\pmod{2^{v_2(2m)}}.
    \end{equation}
    Therefore,
    \begin{align*}
        X_m&\equiv-\sum_{k=0}^{m-1}\binom{2m}{2k}\frac{H^{(2m-2k)}}{H}X_k\quad(\text{remark; }X_0=0)\\
        &\equiv
        \begin{cases}
            \displaystyle2^{v_2(m)}(1+T^2)-\sum_{\substack{k=1\\k\text{:even}}}^{m}\binom{2m}{2k}2^{v_2(k)}(1+T^2)-T\sum_{\substack{k=0\\k\text{:odd}}}^{m}\binom{2m}{2k}(1+T)&(\text{if $m\ge3$ : even})\\
            \displaystyle1+T-\sum_{\substack{k=1\\k\text{:even}}}^{m}\binom{2m}{2k}2^{v_2(k)}T(1+T^2)-\sum_{\substack{k=0\\k\text{:odd}}}^{m}\binom{2m}{2k}(1+T)&(\text{if $m\ge3$ : odd})\\
        \end{cases}\\
        &\equiv
        \begin{cases}
            \displaystyle2^{v_2(m)}(1+T^2)-\sum_{\substack{k=1\\k\text{:even}}}^{m}\binom{2m}{2k}2^{v_2(k)}(1+T^2)&(\text{if $m\ge3$ : even})\\
            \displaystyle1+T\quad(\text{since }2^{v_2(2m)}=2)&(\text{if $m\ge3$ : odd})\\
        \end{cases}\mod{2^{v_2(2m)}}\\
    \end{align*}
    where we used $\sum_{k\text{:odd}}\binom{2m}{2k}\equiv 0\mod{2^{v_2(2m)}}$ for $m\ge 3$ in the last equivalence.
    It is sufficient to prove the following congruence;
    \begin{equation}
    \label{Congruence_for_sum_with_2power}
        \sum_{l=1}^{m/2}\binom{2m}{4l}2^{v_2(2l)}\equiv0\mod{2^{v_2(2m)}}\text{ (for any even number $m\ge 4$)}.
    \end{equation}
    
    In fact,
    \begin{align*}
        2\sum_{l=1}^{m/2}\binom{2m}{4l}2^{v_2(2l)}&=2\sum_{l=1}^{m/2}\binom{2m-1}{4l-1}\frac{2m}{4l}2^{v_2(2l)}\\
        &=2m\sum_{l=1}^{m/2}\binom{2m-1}{4l-1}\frac{2^{v_2(l)}}{l}\\
        &=2m\sum_{l=1}^{m/2}\binom{2m-1}{4l-1}\frac{1}{l'}\text{ (where $l'=l/2^{v_2(l)}$ is odd)}\\
        &\equiv 2m\sum_{l=1}^{m/2}\binom{2m-1}{4l-1} \pmod{2^{v_2(2m)+1}}\\
        &=2m\sum_{l=0}^{2m-1}\binom{2m-1}{l}\frac{1}{4}\sum_{j=0}^{3}(\sqrt{-1})^{j(l+1)}\\
        &=\frac{2m}{4}\sum_{j=0}^{3}(\sqrt{-1})^j(1+(\sqrt{-1})^j)^{2m-1}\\
        &=\frac{2m}{4}\left(2^{2m-1}+\frac{\sqrt{-1}}{1+\sqrt{-1}}(2\sqrt{-1})^m+\frac{-\sqrt{-1}}{1-\sqrt{-1}}(-2\sqrt{-1})^m\right)\\
        &=\frac{2m}{4}(2^{2m-1}+(2\sqrt{-1})^m)\text{ (Note that $m$ is even)}\\
        &=2m(2^{2m-3}-(2\sqrt{-1})^{m-2})\equiv0\mod{2^{v_2(2m)+1}}\text{ (since $m\ge4$)}.
    \end{align*}
    This means (\ref{Congruence_for_sum_with_2power}) holds and proof for the case of $p=2$ is completed.
    
    \underline{In the case of $p=3$}
    
    The main points of the argument are almost the same as in the case of $p=2$. Consider the case $N=3m$ for some $m>0$ in (\ref{Equ_not _congru}). Again, if $k$ isn't a multiple of $3$, $\binom{3m}{k}\equiv 0\mod{3^{v_3(3m)}}$ by Lemma \ref{Lem_binom_p_valuation}. So we have
    \begin{align*}
        \sum_{k=0}^{3m}\binom{3m}{3k}\frac{H^{(3m-3k)}}{H}X_k&=\sum_{k=0}^{3m}\binom{3m}{3k}(-1)^{3k}\frac{H^{(3m-3k)}}{H}\left(1-D_{3k}\right)\quad(\text{definition of }X_k)\\
        &\equiv\sum_{k=0}^{3m}\binom{3m}{k}(-1)^{k}\frac{H^{(3m-k)}}{H}\left(1-D_k\right)\mod{3^{v_3(3m)}}\\
        &=\sum_{k=0}^{3m}\binom{3m}{k}(-1)^k\frac{H^{(3m-k)}}{H}\quad(\text{by (\ref{Equ_not _congru})})\\
        &\equiv\sum_{k=0}^{3m}\binom{3m}{3k}(-1)^k\frac{H^{(3m-3k)}}{H}\\
        &\equiv
        \begin{cases}
            \displaystyle\sum_{\substack{k=0\\k\text{:even}}}^{m}\binom{3m}{3k}-\frac{H^{(3)}}{H}\sum_{\substack{k=0\\k\text{:odd}}}^{m}\binom{3m}{3k}\quad \text{(if }m>0\text{ : even)}\\
            \displaystyle\frac{H^{(3)}}{H}\sum_{\substack{k=0\\k\text{:even}}}^{m}\binom{3m}{3k}-\sum_{\substack{k=0\\k\text{:odd}}}^{m}\binom{3m}{3k}\quad \text{(if }m>0\text{ : odd)}\\
        \end{cases}\mod{3^{v_3(3m)}}.\\
    \end{align*}
    
    Using $\binom{3m}{k}\equiv 0\mod{3^{v_3(3m)}}$ for $3\nmid k$ here too, we have 
    \[\sum_{\substack{k=0\\k\text{:even}}}^{m}\binom{3m}{3k}\equiv\sum_{\substack{k=0\\k\text{:even}}}^{3m}\binom{3m}{k}=2^{3m-1}=\sum_{\substack{k=0\\k\text{:odd}}}^{3m}\binom{3m}{k}\equiv\sum_{\substack{k=0\\k\text{:odd}}}^{m}\binom{3m}{3k}\mod{3^{v_3(3m)}}.\]
    We obtain the following recurrence equation in the end;
    \begin{equation}
    \label{Eq_reccurence}
    \sum_{k=0}^{m}\binom{3m}{3k}\frac{H^{(3m-3k)}}{H}X_m\equiv(-1)^{m}\left(1-\frac{H^{(3)}}{H}\right)2^{3m-1}\equiv\left(T-1\right)(-2)^{3m-1}\mod{3^{v_3(3m)}}.
    \end{equation}
    
    We will show (\ref{Equ_Xm_3}) by induction.
    If $m=1$, $v_3(3m)=1$ and the congruence 
    \[\binom{3}{0}\frac{H^{(3)}}{H}X_0+\binom{3}{3}X_1\equiv \left(T-1\right)(-2)^{3-1}\mod{3}\]
    holds. Since $X_0=0$, we have $X_1\equiv T-1\mod{3}$. 
    On the other hand,
    \[\frac{(1-3T)(T-1)}{1+3T^2}=(1-3T)(T-1)(1-3T^2+9T^4-\cdots)\equiv1\cdot (T-1)\cdot1\mod{3}.\]
    So (\ref{Equ_Xm_3}) holds for $m=1$.

    Assume that (\ref{Equ_Xm_3}) holds for $1\le k< m$. Then, by Lemma \ref{Lem_binom_p_valuation}, 
    \begin{equation}
    \label{Eq_asuump}
        \binom{3m}{3k}X_k\equiv
        \begin{cases}
            \displaystyle\binom{3m}{3k}\frac{(1+T)(T-1)}{1+3T^2} &(\text{if }k\ge1\text{ : even})\\
            \displaystyle\binom{3m}{3k}\frac{(1-3T)(T-1)}{1+3T^2} &(\text{if }k\ge1\text{ : odd})\\
        \end{cases}
        \quad\pmod{3^{v_3(3m)}}
    \end{equation}
    for all $1\le k<m$.
    If $m$ is even, substituting (\ref{Eq_asuump}) into (\ref{Eq_reccurence}) yields
    \begin{align*}
        (-2)^{3m-1}(T-1)&\equiv\sum_{k=0}^{m}\binom{3m}{3k}\frac{H^{(3m-3k)}}{H}X_k\\
        &\equiv\sum_{\substack{k=0\\k\text{:odd}}}^m\binom{3m}{3k}T\frac{(1-3T)(T-1)}{1+3T^2}+\sum_{\substack{k=0\\k\text{:even}}}^m\binom{3m}{3k}\frac{(1+T)(T-1)}{1+3T^2}-2\frac{(1+T)(T-1)}{1+3T^2}+X_m\\
        &\equiv2^{3m-1}\left(T\frac{(1-3T)(T-1)}{1+3T^2}+\frac{(1+T)(T-1)}{1+3T^2}\right)-2\frac{(1+T)(T-1)}{1+3T^2}+X_m\\
        &=2^{3m-1}\frac{(1+2T-3T^2)(T-1)}{1+3T^2}-2\frac{(1+T)(T-1)}{1+3T^2}+X_m\\
        &=2^{3m-1}\left(2\frac{(1+T)(T-1)}{1+3T^2}-(T-1)\right)-2\frac{(1+T)(T-1)}{1+3T^2}+X_m\\
        &=(3-1)^{3m}\frac{(1+T)(T-1)}{1+3T^2}+(-2)^{3m-1}(T-1)-2\frac{(1+T)(T-1)}{1+3T^2}+X_m\\
        &(-2^{3m-1}=(-2)^{3m-1}\text{ since }m\text{ : even})\\
        &\equiv(-1)^{3m}\frac{(1+T)(T-1)}{1+3T^2}+(-2)^{3m-1}(T-1)-2\frac{(1+T)(T-1)}{1+3T^2}+X_m\\
        &=(-2)^{3m-1}(T-1)+X_m-\frac{(1+T)(T-1)}{1+3T^2}\quad\mod{3^{v_3(3m)}}.
    \end{align*}
    Thus, we have $X_m\equiv\frac{(1+T)(T-1)}{1+3T^2}\mod{3^{v_3(3m)}}$.

    If $m$ is odd, similarly we have
    \begin{align*}
        (-2)^{3m-1}(T-1)&\equiv\sum_{k=0}^{m}\binom{3m}{3k}\frac{H^{(3m-3k)}}{H}X_k\\
        &\equiv\sum_{\substack{k=0\\k\text{:odd}}}^m\binom{3m}{3k}\frac{(1-3T)(T-1)}{1+3T^2}+\sum_{\substack{k=0\\k\text{:even}}}^m\binom{3m}{3k}T\frac{(1+T)(T-1)}{1+3T^2}\\
        &\qquad-T\frac{(1+T)(T-1)}{1+3T^2}-\frac{(1-3T)(T-1)}{1+3T^2}+X_m\\
        &\equiv2^{3m-1}\frac{(1-2T+T^2)(T-1)}{1+3T^2}-T\frac{(1+T)(T-1)}{1+3T^2}-\frac{(1-3T)(T-1)}{1+3T^2}+X_m\\
        &\equiv2^{3m-1}((T-1)-2T\frac{(1+T)(T-1)}{1+3T^2})-T\frac{(1+T)(T-1)}{1+3T^2}-\frac{(1-3T)(T-1)}{1+3T^2}+X_m\\
        &\equiv(-2)^{3m-1}(T-1)+((3-1)^{3m}-1)T\frac{(1+T)(T-1)}{1+3T^2}+X_m-\frac{(1-3T)(T-1)}{1+3T^2}\\
        &(2^{3m-1}=(-2)^{3m-1}\text{ since }m\text{ : odd})\\
        &\equiv(-2)^{3m-1}(T-1)+X_m-\frac{(1-3T)(T-1)}{1+3T^2}\quad\mod{3^{v_3(3m)}}
    \end{align*}
    Hence, (\ref{Equ_Xm_3}) also holds for $m$.
\end{proof}

\begin{lemma}
\label{Lem_2_for_special_case}
    Let $H(z)=\sum_{n=0}^{\infty}z^{6n}/(6n)!$. 
    
    (1) For all $r>0$, there exists some polynomial $P_r(z)$ such that
    \[(H^{(3)})^2-H^2\equiv P_r(z)\mod{3^r}.\]

    (2) For all $r>0$, there exists some polynomial $Q_r(z)$ such that
    \[\frac{1}{H^3+3H(H^{(3)})^2}\equiv Q_r\mod{3^r}.\]
\end{lemma}
\begin{proof}
    (1) We have 
    \begin{align*}
        (H^{(3)}(z))^2-H(z)^2&=\left(\sum_{n=0}^{\infty}\frac{z^{6n+3}}{(6n+3)!}\right)^2-\left(\sum_{n=0}^{\infty}\frac{z^{6n}}{(6n)!}\right)^2\\
        &=\sum_{n=1}^{\infty}\left(\sum_{m=0}^{n-1}\binom{6n}{6m+3}\right)\frac{z^{6n}}{(6n)!}-\sum_{n=0}^{\infty}\left(\sum_{m=0}^{n}\binom{6n}{6m}\right)\frac{z^{6n}}{(6n)!}\\
        &=-1-\sum_{n=1}^{\infty}\left(\sum_{k=0}^{2n}(-1)^k\binom{6n}{3k}\right)\frac{z^{6n}}{(6n)!}\\
        &=-1-\sum_{n=1}^{\infty}\left(\sum_{k=0}^{2n}(-1)^k\binom{6n}{k}\frac{1}{3}\sum_{j=0}^{2}\omega^{jk}\right)\frac{z^{6n}}{(6n)!}\quad(\omega\coloneqq\frac{-1+\sqrt{-3}}{2})\\
        &=-1-\sum_{n=1}^{\infty}\left(\frac{1}{3}\sum_{j=0}^{2}(1-\omega^j)^{6n}\right)\frac{z^{6n}}{(6n)!}\\
        &=-1-\sum_{n=1}^{\infty}\left(\frac{1}{3}(0^{6n}+(\sqrt{-3}\omega^2)^{6n}+(\sqrt{-3}\omega)^{6n})\right)\frac{z^{6n}}{(6n)!}\quad(\omega\coloneqq\frac{-1+\sqrt{-3}}{2})\\
        &=-1-\sum_{n=1}^{\infty}\left((-1)^n\cdot2\cdot3^{3n-1}\right)\frac{z^{6n}}{(6n)!}\\
        &\equiv P_r\mod{3^r}
    \end{align*}
    where $P_r$ is a polynomial depending on $r$.

    (2) First, we identify the $3$-adic valuation of the coefficients $\{c_{n}\}$ of $H^3+3H(H^{(3)})^2=\sum_{n=0}^{\infty}c_nz^{6n}/(6n)!$. 
    
    We use $G(z)\coloneqq H_{3,0}(z)=\sum_{n=0}^{\infty}z^{3n}/(3n)!$ here. In terms of $G$,
    \begin{align*}
        H(z)^3+3H(z)(H^{(3)}(z))^2&=H(z)(H(z)^2+3H^{(3)}(z)^2)\\
        &=\left(\frac{G(z)+G(-z)}{2}\right)\left(\left(\frac{G(z)+G(-z)}{2}\right)^2+3\left(\frac{G(z)-G(-z)}{2}\right)^2\right)\\
        &=\left(\frac{G(z)+G(-z)}{8}\right)(4G(z)^2-4G(z)G(-z)+4G(-z)^2)\\
        &=\frac{G(z)^3+G(-z)^3}{2}.
    \end{align*}
    Here, it becomes simpler to express $G$ as a sum of exponents.
    \begin{align*}
        G(z)^3&=\left(\sum_{n=0}^{\infty}\frac{z^{3n}}{(3n)!}\right)^3\\
        &=\frac{1}{27}(e^z+e^{\omega z}+e^{\omega^2z})^3\quad(\omega=\frac{-1+\sqrt{-3}}{2})\\
        &=\frac{1}{27}(e^{3z}+e^{3\omega z}+e^{3\omega^2z}+3e^{(1+2\omega)z}+3e^{(2+\omega)z}\\
        &+3e^{(1+2\omega^2)z}+3e^{(2+\omega^2)z}+3e^{(\omega+2\omega^2)z}+3e^{(2\omega+\omega^2)z}+6).
    \end{align*}
    Since 
    $\begin{cases}
        -(1+2\omega^2)=1+2\omega=\sqrt{-3}\\
            -(2+\omega)=\omega+2\omega^2=\sqrt{-3}\omega\\
            -(2\omega+\omega^2)=2+\omega^2=\sqrt{-3}\omega^2
    \end{cases},$
    we have
    \begin{align*}
        \frac{G(z)^3+G(-z)^3}{2}&=\frac{G(3z)+G(-3z)}{18}+\frac{G(\sqrt{-3}z)+G(-\sqrt{-3}z)}{3}+\frac{2}{9}\\
        &=1+\sum_{n=1}^{\infty}\left(\frac{2}{9}3^{6n}+\frac{2}{3}(\sqrt{-3})^{6n}\right)\frac{z^{6n}}{(6n)!}
    \end{align*}
    Thus, $\{c_n\}$, the coefficients of $H(z)^3+3H(z)(H^{(3)}(z))^2$ satisfy $v_3(c_n)=3n-1$.

    Let ${d_n}_{n=0}^{\infty}$ be defined by
    \[
    \sum_{n=0}^{\infty} d_n \frac{z^n}{n!}=\frac{1}{H(z)^3+3H(z)(H^{(3)}(z))^2}.
    \]
    Then, we have 
    \[\sum_{m=0}^{n}\binom{6n}{6m}d_mc_{n-m}=
    \begin{cases}
        1&(n=0)\\
        0&(n>0)
    \end{cases}
    \]
    by definition. Because $c_0=1$,
    \[v_3(d_n)=v_3\left(-\sum_{m=0}^{n-1}\binom{6n}{6m}d_mc_{n-m}\right)\ge\min_{0\le m<n}(v_3(d_m)+v_3(c_{n-m}))\quad(\text{for all }n>0).\]
    Then, we can conclude $v_3(d_n)\ge 2n$ by induction. This implies $1/(H^3+3H(H^{(3)})^2)\equiv Q_r\mod{3^r}$ for some polynomial $Q_r$.
\end{proof}

Note that the above proof for Proposition \ref{Prop_for_conj_of_special_case} cannot be generalized to the case of $j>0$ immediately. In fact, the nonvanishing of the constant term plays an important role. If this is not the case, we have Laurent series and we need to take the term $1/z^M$ into account in the end. However, this causes a shift in the coefficients and this is crucial because we work with $\sum (a_n/n!)z^n$ rather than $\sum a_nz^n$.

\section{Application: Even Zeta values and Bernoulli numbers}

As mentioned in the introduction, in the case $N=2$ and $j=0$, $\mathcal{E}_{Nn}^{(N,j)}$ are Euler numbers and are thus related to zeta values and Bernoulli numbers. Recall that there is a well known formula for Euler numbers and Bernoulli numbers:
\[E_{2n}=\sum_{k=1}^{n}\binom{2n}{2k-1}\frac{2^{2k}-4^{2k}}{2k}B_{2k}+1\]
\[B_{2n}=\frac{2n}{4^{2n}-2^{2n}}\sum_{k=0}^{n-1}\binom{2n-1}{2k}E_{2k}.\]

Similarly, in the cases of $(N,j)=(4,0),(4,2),(6,3)$, the numbers $\{\mathcal{E}_{Nn}^{(N,j)}\}$ are also related to zeta values and hence to Bernoulli numbers. Which cases of congruential Euler numbers $\{\mathcal{E}_{Nn}^{(N,j)}\}$ admit such connections is determined by the distribution of the zeros of $H_{N,j}(z)=\sum_{n=0}^{\infty}z^{Nn+j}/(Nn+j)!=N^{-1}\sum_{k=0}^{N-1}\zeta_{N}^{-kj}\exp{(\zeta_N^{k}z)}$.

In this section, using complex analysis, we describe how congruential Euler numbers represent even zeta values.

The following proposition shows that the zeros of $H_{N,j}(z)$ are regularly distributed in the cases of $(N,j)=(4,0),(4,2), \text{and }(6,3)$.
\begin{proposition}
\label{Prop_zeros_of_H}
    The nontrivial zeros (i.e. zeros except for $z=0$) of $H_{N,j}(z)$ are
    \begin{align*}
        &z_{k,l}^{(4,0)}=\sqrt{2}\zeta_8\zeta_4^l(k+\frac{1}{2})\pi\quad(k>0,\le l<4)\quad\text{if }(N,j)=(4,0)\\
        &z_{k,l}^{(4,2)}=\sqrt{2}\zeta_8\zeta_4^lk\pi\quad(k>0,\le l<4)\quad\text{if }(N,j)=(4,2)\\
        &z_{k,l}^{(6,3)}=2\zeta_{12}\zeta_6^lk\pi\quad(k>0,\le l<6)\quad\text{if }(N,j)=(6,3).\\
    \end{align*}
\end{proposition}
\begin{proof}
    The statement follows from the following expressions of $H_{N,j}$ by trigonometric functions.
    \begin{align*}
        H_{4,0}(z)&=\frac{1}{4}(e^z+e^{\sqrt{-1}z}+e^{-z}+e^{-\sqrt{-1}z})\\
        &=\frac{1}{2}(\cos{z}+\cos{(\sqrt{-1}z)})\\
        &=\cos{(\frac{1+\sqrt{-1}}{2}z)}\cos{(\frac{1-\sqrt{-1}}{2}z)}
    \end{align*}
    \begin{align*}
        H_{4,2}(z)&=\frac{1}{4}(e^z-e^{\sqrt{-1}z}+e^{-z}-e^{-\sqrt{-1}z})\\
        &=\frac{1}{2}(-\cos{z}+\cos{(\sqrt{-1}z)})\\
        &=\sin{(\frac{1+\sqrt{-1}}{2}z)}\sin{(\frac{1-\sqrt{-1}}{2}z)}
    \end{align*}
    \begin{align*}
        H_{6,3}(z)&=\frac{1}{6}(e^z-e^{\zeta_6z}+e^{\zeta_6^2z}-e^{\zeta_6^3z}+e^{\zeta_6^4z}-e^{\zeta_6^5z})\\
        &=-\frac{\sqrt{-1}}{3}(\sin{(\sqrt{-1}z)}-\sin{(\zeta_6\sqrt{-1}z)}+\sin{(\zeta_6^2\sqrt{-1}z)})\\
        &=-\frac{2\sqrt{-1}}{3}\left(\sin{(\frac{\sqrt{-1}}{2}z)}\cos{(\frac{\sqrt{-1}}{2}z)}+\sin{(\frac{\zeta_6^2-\zeta_6}{2}\sqrt{-1}z)}\cos{(\frac{\zeta_6+\zeta_6^2}{2}\sqrt{-1}z)}\right)\\
        &=-\frac{2\sqrt{-1}}{3}\left(\sin{(\frac{\sqrt{-1}}{2}z)}\cos{(\frac{\sqrt{-1}}{2}z)}-\sin{(\frac{\sqrt{-1}}{2}z)}\cos{(\frac{\sqrt{3}}{2}z)}\right)\\
        &=-\frac{2\sqrt{-1}}{3}\sin{(\frac{\sqrt{-1}}{2}z)}\sin{(\frac{\sqrt{3}+\sqrt{-1}}{4}z)}\sin{(\frac{\sqrt{3}-\sqrt{-1}}{4}z)}
    \end{align*}
\end{proof}

\begin{remark}
    For general $N$, the nontrivial zeros of $H_{N,j}(z)$ are expected to be approximately located on the lines $\arg{z}=\frac{\pi}{N}+\frac{2\pi k}{N}$, but strict regularity like Proposition \ref{Prop_zeros_of_H} is not observed for the other pairs $(N,j)$. One might expect that the zeros of $H_{2p,p}$ have strict regularity, however, computations show otherwise. For example, computer calculations suggest
    \[\left|\frac{\mathcal{E}_{10n}^{(10,5)}}{(10n)!}\middle/\frac{\mathcal{E}_{10n+10}^{(10,5)}}{(10n+10)!}\right|/\pi^{10}\rightarrow116471.365616430\sim(3.210864043359495)^{10}\quad(n\rightarrow\infty).\]
     By d’Alembert’s ratio test, this means the distance from the origin to the closest nontrivial zero of $H_{10,5}$ is about $(3.210864)\pi$. The value $3.210864\cdots$ does not seem to have a simple closed form; in particular, it is unclear whether it is algebraic.
\end{remark}

In the representation of even zeta values using congruent Euler numbers, the following lemma allows us to obtain the expression for $\zeta(k)$ for even $k$ that are not multiples of $6$.

\begin{lemma}
\label{Lem_special_values_of_H62_H64}
    (1) For any $k>0$ and $0\le l<4$,
    \begin{equation}
    \label{Equ_special_values_for_(40)}
        H_{4,1}(z_{k,l}^{(4,0)})=\zeta_4^{2l+3}H_{4,3}(z_{k,l}^{(4,0)}),
    \end{equation}
    \begin{equation}
    \label{Equ_special_values_for_(42)}
        H_{4,3}(z_{k,l}^{(4,2)})=\zeta_4^{2l+3}H_{4,1}(z_{k,l}^{(4,2)}).
    \end{equation}
    
    (2) For any $k>0$ and $0\le l<6$,
    \begin{equation}
    \label{Equ_special_values_for_(63)}
        H_{6,4}(z_{k,l}^{(6,3)})=\zeta_6^{2(l-1)}H_{6,2}(z_{k,l}^{(6,3)}).
    \end{equation}
\end{lemma}
\begin{proof}
    In general, 
    \begin{equation*}
        H_{N,j}(\zeta_Nz)=\frac{1}{N}\sum_{k=0}^{N-1}\zeta_N^{-kj}e^{\zeta_N^{k+1}z}=\frac{1}{N}\sum_{k=0}^{N-1}\zeta_N^{-(k-1)j}e^{\zeta_N^{k}z}=\zeta_N^jH_{N,j}(z).
    \end{equation*}
    Therefore, for any $z_0\in\mathbb{C}$,
    \begin{equation}
    \label{Eq_fixed_l0}
    \begin{aligned}
        &\text{if}\quad H_{N,j}(\zeta_N^{l_0}z_0)=cH_{N,j'}(\zeta_N^{l_0}z_0) \quad\text{for some }l_0\text{ and }c,\\
        &\text{then}\quad H_{N,j}(\zeta_N^lz_0)=c\zeta_N^{(j-j')(l-l_0)}H_{N,j'}(\zeta_N^lz_0)\quad\text{for all }l.
    \end{aligned}
    \end{equation}
    
    (1) Substituting $z=z_{k,l}^{(4,0)}=(1+\sqrt{-1})(k-\frac{1}{2})\pi$ and $z=z_{k,l}^{(4,2)}=(1+\sqrt{-1})k\pi$ for $H_{4,1}(z)$ and $H_{4,3}(z)$, we have
    \begin{align*}
        &H_{4,1}((1+\sqrt{-1})(k-\frac{1}{2})\pi)=(-1)^{k-1}(1+\sqrt{-1})(e^{(k-\frac{1}{2})\pi}+e^{-(k-\frac{1}{2})\pi}),\\
        &H_{4,3}((1+\sqrt{-1})(k-\frac{1}{2})\pi)=(-1)^{k-1}(-1+\sqrt{-1})(e^{(k-\frac{1}{2})\pi}+e^{-(k-\frac{1}{2})\pi}),\\
        &H_{4,1}((1+\sqrt{-1})k\pi)=(-1)^k(1+\sqrt{-1})(e^{k\pi}-e^{-k\pi}),\\
        &H_{4,3}((1+\sqrt{-1})k\pi)=(-1)^k(1-\sqrt{-1})(e^{k\pi}-e^{-k\pi}).
    \end{align*}
    Thus, we obtain $H_{4,1}(z_{k,0}^{(4,0)})=-\sqrt{-1}H_{4,3}(z_{k,0}^{(4,0)})$ and $H_{4,3}(z_{k,0}^{(4,2)})=-\sqrt{-1}H_{4,1}(z_{k,0}^{(4,2)})$. These equations lead to \eqref{Equ_special_values_for_(40)} and \eqref{Equ_special_values_for_(42)} by the claim \eqref{Eq_fixed_l0}.
    
    (2) Let $j=2$ or $4$. Then, we have
    \begin{align*}
        H_{6,j}(z)&=\frac{1}{6}\sum_{k=0}^{5}\zeta_6^{-kj}e^{\zeta_6^kz}\\
        &=\frac{1}{6}\sum_{k=0}^{5}\omega^{-kj/2}e^{\zeta_6^kz}\quad(\omega=\zeta_6^2=\frac{-1+\sqrt{-3}}{2})\\
        &=\frac{1}{3}(\cos{(\sqrt{-1}z)}+\omega^{-j/2}\cos{(\zeta_6\sqrt{-1}z)}+\omega^{-j}\cos{(\zeta_6^2\sqrt{-1}z)})\\
        &=\frac{1}{3}\left(\cos{(\sqrt{-1}z)}-\frac{1}{2}(\cos{(\zeta_6\sqrt{-1}z)}+\cos{(\zeta_6^2\sqrt{-1}z)})\right.\\
        &\qquad\qquad\qquad\left.+(-1)^{j/2}\frac{\sqrt{-3}}{2}(\cos{(\zeta_6\sqrt{-1}z)}-\cos{(\zeta_6^2\sqrt{-1}z)})\right)\\
        &=\frac{1}{3}\left(\cos{(\sqrt{-1}z)}-\cos{(\frac{\zeta_6+\zeta_6^2}{2}\sqrt{-1}z)}\cos{(\frac{\zeta_6-\zeta_6^2}{2}\sqrt{-1}z)}\right.\\
        &\qquad\qquad\qquad\left.-(-1)^{j/2}\sqrt{-3}\sin{(\frac{\zeta_6+\zeta_6^2}{2}\sqrt{-1}z)}\sin{(\frac{\zeta_6-\zeta_6^2}{2}\sqrt{-1}z)}\right)\\
        &=\frac{1}{3}\left(\cos{(\sqrt{-1}z)}-\cos{(\frac{\sqrt{-1}}{2}z)}\cos{(\frac{\sqrt{3}}{2}z)}+(-1)^{j/2}\sin{(\frac{\sqrt{-1}}{2}z)}\sin{(\frac{\sqrt{3}}{2}z)}\right).
    \end{align*}
    Substituting $z=z_{k,1}^{(6,3)}=2k\pi\sqrt{-1}$ for this, we obtain
    \[H_{6,j}(z_{k,1}^{(6,3)})=\frac{1-(-1)^k\cosh{(\sqrt{3}k\pi)}}{3}\quad(j=2,4).\]
    
    So we have \eqref{Equ_special_values_for_(63)} by \eqref{Eq_fixed_l0}.
\end{proof}

\begin{lemma}
    Let $(N,j)=(4,0),\text{ }(4,2),\text{ or }(6,3)$. For any integers $0\le j'<N$ and $M\ge2$ ,
    \[\int_{C_n}\frac{H_{N,j'}(z)}{z^MH_{N,j}(z)}dz\rightarrow0\quad\text{when }n\rightarrow\infty\]
    where $C_n$ is the circle $|z|=R_n\text{ }(n\in\mathbb{Z}_{>0})$ and $R_n=\begin{cases}
        \sqrt{2}n\pi&\text{if }(N,j)=(4,0)\\
        \sqrt{2}(n-\frac{1}{2})\pi&\text{if }(N,j)=(4,2)\\
        2(n-\frac{1}{2})\pi&\text{if }(N,j)=(6,3)\\
    \end{cases}.$
\end{lemma}
\begin{proof}
    It suffices to show that there exists a constant $C$ such that $|H_{N,j'}(z)/H_{N,j}(z)|<C$ for any $z$ on $C_n$ of any sufficiently large $n$. In fact, if it is true, we have
    \[\left|\int_{C_n}\frac{H_{N,j'}(z)}{z^MH_{N,j}(z)}dz\right|\le\int_{C_n} \frac{C}{R_n^M}|dz|=\frac{C}{R_n^{M-1}}\rightarrow0\text{ }(n\rightarrow\infty).\]
    
    Roughly speaking, $|H_{N,J}(z)|\text{ }(0\le J<N)$ is dominated by $|e^{\zeta_N^Kz}|$ where $K$ satisfies \[\Re({\zeta_N^Kz)}=|z|\cos{(\arg{z}+\frac{2\pi K}{N})}=\max_{k}{\Re{(\zeta_N^k}z)}=|z|\max_k\cos{(\arg{z}+\frac{2\pi k}{N})}\] where $\Re$ denotes the real part. This $K$ does not depend on $J$, so $|H_{N,j'}(z)/H_{N,j}(z)|$ is bounded except for around zeros of $H_{N,j}$.
    
    We will give a rigorous proof.
    
    Since $H_{N,J}(\zeta_Nz)=\zeta_N^JH(z)\text{ and }H_{N,J}(\overline{z})=\overline{H_{N,J}(z)}\text{ }(0\le J<N)$, it suffices to show that there are some constants $C>0$ such that \[|H_{N,j'}(z)/H_{N,j}(z)|<C\] for any sufficiently large $n$ and any $z\in C_n$ satisfying $-\frac{\pi}{N}\le\theta\coloneqq\arg(z)\le0$.

    For $z$ with $-\frac{\pi}{N}\le\theta\le0$, the main term of $H_{N,j}(z)$ is $\frac{1}{N}e^z$ and the second largest term is $\frac{1}{N}e^{\zeta_Nz}$. 
    
    Let 
    \[\varepsilon=\begin{cases}
        1 & (\text{if }(N,j)=(4,0))\\
        -1 & (\text{if }(N,j)=(4,2)\text{ or }(6,3))
    \end{cases}.\]
    Then, $H_{N,J}\text{ }(0\le J<N)$ have the form 
    \[H_{N,J}(z)=\frac{1}{N}(e^z+\varepsilon e^{\zeta_Nz}+(\text{the other terms}))\]
    and $e^z+\varepsilon e^{\zeta_Nz}$ is the main term.

    First, we show that this main term does not vanish on $\{z\mid z\in C_n\text{ for some }n\text{ with }-\frac{\pi}{N}\le\theta\le0\}$. We have
    \[|e^z+\varepsilon e^{\zeta_Nz}|=|e^z||1+\varepsilon e^{x+\sqrt{-1}y}|\quad\text{where }x+\sqrt{-1}y\coloneqq(\zeta_N-1)z\quad(x,y\in\mathbb{R}).\]
    Because
    \[x^2+y^2=|\zeta_N-1|^2|z|^2=
    \begin{cases}
        (2n\pi)^2 & (\text{if }(N,j)=(4,0))\\
        ((2n-1)\pi)^2 & (\text{if }(N,j)=(4,2)\text{ or }(6,3))
    \end{cases},\] 
    there exists some $\delta>0$ such that for any $z\in C_n$ with $|x|<\delta$,
    \[\begin{cases}
        |e^{x+\sqrt{-1}y}-1|<1 & (\text{if }(N,j)=(4,0))\\
        |e^{x+\sqrt{-1}y}-(-1)|<1 & (\text{if }(N,j)=(4,2)\text{ or }(6,3))
    \end{cases}
    \quad\text{i.e. }|e^{x+\sqrt{-1}y}-\varepsilon|<1.\]
    In particular,
    \[|1+\varepsilon e^{x+\sqrt{-1}y}|=|e^{x+\sqrt{-1}y}-(-\varepsilon)|>1\quad\text{for all }|x|<\delta.\]
    On the other hand, since now we consider $z$ with $-\frac{\pi}{N}\le\theta\le0$, 
    \[\begin{cases}
        \frac{\pi}{2}\le\arg{(x+\sqrt{-1}y)}\le\frac{3\pi}{4} & (\text{if }N=4)\\
        \frac{\pi}{2}\le\arg{(x+\sqrt{-1}y)}\le\frac{2\pi}{3} & (\text{if }N=6)
    \end{cases}.\]
    Hence for $z$ with $|x|\ge\delta$, we have $x<-\delta$.
    Thus, 
    \[|1+\varepsilon e^{x+\sqrt{-1}y}|\ge|1-|e^{x+\sqrt{-1}y}||=|1-e^x|\ge1-e^{-\delta}.\]
    Therefore there exists a constant $C_0>0$ such that 
    \[|e^z+\varepsilon e^{\zeta_Nz}|>C_0|e^z|\quad\text{for all }z\text{ with }-\frac{\pi}{N}<\theta\le0.\]
    In particular, $e^z$ and $e^{\zeta_Nz}$ do not cancel each other, and therefore there exists a constant $C_1>0$ and a positive integer $n_0$ such that 
    \[|H_{N,j}(z)|>C_1|e^z|\quad\text{for all }z\text{ on }C_n(n>n_0)\text{ with }-\frac{\pi}{N}\le\theta\le0.\]

    On the other hand, since $|e^z|$ is equal to or larger than $|e^{\zeta_kz}|$ for all $0<k<N$ when $-\frac{\pi}{N}\le\theta\le0$, there exists a constant $C_2>0$ and a positive integer $n_1$ such that 
    \[|H_{N,j'}(z)|<C_2|e^z|\quad\text{for all }z\text{ on }C_n(n>n_1)\text{ with }-\frac{\pi}{N}\le\theta\le0.\]

    Finally, we can conclude that $|H_{N,j'}(z)/H_{N,j}(z)|$ is bounded for any sufficiently large $n$ and any $z$ on $C_n$ with $-\frac{\pi}{N}\le\theta\le0$.
\end{proof}

\begin{theorem}[Restatement of Theorem \ref{MainThm2}]
    For any positive integer $n$, the following equations hold.
    \begin{equation}
    \label{Equ_E40_and_lambda}
        \lambda(4n)=\sum_{k=1}^\infty\frac{1}{(2k-1)^{4n}}=\frac{(-1)^{n+1}(\pi/\sqrt{2})^{4n}}{4(4n-1)!}\sum_{m=0}^{n-1}\binom{4n-1}{4m}\mathcal{E}_{4m}^{(4,0)},
    \end{equation}
    \begin{equation}
    \label{Equ_E40_and_lambda_2}
        \lambda(4n-2)=\sum_{k=1}^\infty\frac{1}{(2k-1)^{4n-2}}=\frac{(-1)^{n+1}(\pi/\sqrt{2})^{4n-2}}{4(4n-3)!}\sum_{m=0}^{n-1}\binom{4n-3}{4m}\mathcal{E}_{4m}^{(4,0)},
    \end{equation}
    \begin{equation}
    \label{Equ_E40_and_zeta}
        \zeta(4n)=\frac{(-1)^{n+1}(\sqrt{2}\pi)^{4n}}{4(4n-1)!(2^{4n}-1)}\sum_{m=0}^{n-1}\binom{4n-1}{4m}\mathcal{E}_{4m}^{(4,0)},
    \end{equation}
    \begin{equation}
    \label{Equ_E40_and_zeta_2}
        \zeta(4n-2)=\frac{(-1)^{n+1}(\sqrt{2}\pi)^{4n-2}}{4(4n-3)!(2^{4n-2}-1)}\sum_{m=0}^{n-1}\binom{4n-3}{4m}\mathcal{E}_{4m}^{(4,0)},
    \end{equation}
    \begin{equation}
    \label{Equ_E42_and_zeta}
        \zeta(4n)=\frac{(-1)^{n+1}(\sqrt{2}\pi)^{4n}}{4(4n+1)!}\sum_{m=0}^{n}\binom{4n+1}{4m}\mathcal{E}_{4m}^{(4,2)},
    \end{equation}
    \begin{equation}
    \label{Equ_E42_and_zeta_2}
        \zeta(4n-2)=\frac{(-1)^{n+1}(\sqrt{2}\pi)^{4n-2}}{4(4n-1)!}\sum_{m=0}^{n-1}\binom{4n-1}{4m}\mathcal{E}_{4m}^{(4,2)},
    \end{equation}
    \begin{equation}
    \label{Equ_E63_and_zeta}
        \zeta(6n)=\frac{(-1)^{n+1}(2\pi)^{6n}}{6(6n+2)!}\sum_{m=0}^{n}\binom{6n+2}{6m}\mathcal{E}_{6m}^{(6,3)},
    \end{equation}
    \begin{equation}
    \label{Equ_E63_and_zeta_2}
        \zeta(6n-4)=\frac{(-1)^{n+1}(2\pi)^{6n-4}}{6(6n-2)!}\sum_{m=0}^{n-1}\binom{6n-2}{6m}\mathcal{E}_{6m}^{(6,3)},
    \end{equation}
    where $\zeta(s)$ is the Riemann zeta function and
    $\lambda(s)=(1-2^{-s})\zeta(s)=L(s,\chi)$ (where $\chi$ is the trivial character modulo $2$) is the Dirichlet lambda function.
\end{theorem}
\begin{proof}
    Let $(N,j,j')=(4,0,3),\text{ }(4,0,1),\text{ }(4,2,1),\text{ }(4,2,3),\text{ }(6,3,2),\text{ or }(6,3,4)$. 
    Let $M\ge\max\{2,-j+j'+1\}$ be an integer.
    Then, $f(z)\coloneqq H_{N,j'}(z)/z^MH_{N,j}(z)$ has a pole at $z=0$.
    Let $C_n$ be as in the previous Lemma.
    Using the residue theorem for the integral 
    \[\int_{C_n}\frac{H_{6,j'}(z)}{z^MH_{6,3}(z)}dz\]
    and taking the limit of $n\rightarrow\infty$, we have
    \[\operatorname{Res}(f,0)+\sum_{k=1}^{\infty}\sum_{l=0}^{N-1}\operatorname{Res}(f,z_{k,l}^{(N,j)})=0\]
    where $z_{k,l}^{(N,j)}$ is the one defined in Proposition \ref{Prop_zeros_of_H}.

    $\operatorname{Res}(f,0)$ is the coefficient of $z^{M-1}$ of 
    \begin{align*}
        \frac{H_{N,j'}(z)}{H_{N,j}(z)}&=\frac{1}{z^j}F_{N,j}(z)H_{N,j'}(z)\\
        &=\frac{1}{z^j}\left(\sum_{n=0}^{\infty}\mathcal{E}_{Nn}^{(N,j)}\frac{z^{Nn}}{(Nn)!}\right)\left(\sum_{n=0}^{\infty}\frac{z^{Nn+j'}}{(Nn+j')!}\right)\\
        &=\frac{1}{z^j}\sum_{n=0}^{\infty}\left(\sum_{m=0}^{n}\binom{Nn+j'}{Nm}\mathcal{E}_{Nn}^{(N,j)}\right)\frac{z^{Nn+j'}}{(Nn+j')!}.
    \end{align*}
    Thus,
    \[\operatorname{Res}(f,0)=\begin{cases}
        \displaystyle\frac{1}{(Nn+j')!}\sum_{m=0}^{n}\binom{Nn+j'}{Nm}\mathcal{E}_{Nm}^{(N,j)}&(\text{if }M=Nn-j+j'+1\text{ for some }n)\\
        0&(\text{otherwise})
    \end{cases}.\]
    On the other hand, $z_{k,l}^{(N,j)}$ are simple poles of $f$. Hence, writing $z'=z_{k,l}^{(N,j)}$ and $H_{4,-1}(z)\coloneqq H_{4,3}(z)$ for simplicity, we have
    \begin{align*}
        \operatorname{Res}(f,z_{k,l}^{(N,j)})&=\lim_{z\rightarrow z'}(z-z')\frac{H_{N,j'}(z)}{z^MH_{N,j}(z)}\\
        &=\lim_{z\rightarrow z'}\frac{H_{N,j'}(z)+(z-z')H_{N,j'-1}(z)}{Mz^{M-1}H_{N,j}(z)+z^MH_{N,j-1}(z)}\quad(\text{by l'Hôpital's rule})\\
        &=\begin{cases}
            \displaystyle\frac{1}{z'^M}&(\text{if }(N,j,j')=(4,0,3),\text{ }(4,2,1),\text{ or }(6,3,2))\\
            \displaystyle\frac{\zeta_4^{2l+3}}{z'^M}&(\text{if }(N,j,j')=(4,0,1)\text{ or }(4,2,3)\text{ by Lemma \ref{Lem_special_values_of_H62_H64}})\\
            \displaystyle\frac{\zeta_6^{2(l-1)}}{z'^M}&(\text{if }(N,j,j')=(6,3,4)\text{ by Lemma \ref{Lem_special_values_of_H62_H64}})
        \end{cases}\\
        &=\frac{C_{(N,j,j')}\zeta_N^{l(-j+j'+1)}}{z'^M}\quad
        \left(C_{(N,j,j')}=\begin{cases}
            1 & (\text{if }(N,j,j')=(4,0,3),\text{ }(4,2,1),\text{ or }(6,3,2))\\
            \zeta_4^3 & (\text{if }(N,j,j')=(4,0,1)\text{ or }(4,2,3))\\
            \zeta_6^4 & (\text{if }(N,j,j')=(6,3,4))\\
        \end{cases}\right).
    \end{align*}
    Therefore, if $M=Nn-j+j'+1$ for some $n$, then we  obtain
    \begin{align*}
        \frac{1}{Nn+j'}\sum_{m=0}^{n}\binom{Nn+j'}{Nm}\mathcal{E}_{Nm}^{(N,j)}&=-\sum_{k=1}^{\infty}\sum_{l=0}^{N-1}\operatorname{Res}(f,z_{k,l}^{(N,j)})\\
        &=-\sum_{k=1}^{\infty}\sum_{l=0}^{N-1}\frac{C_{(N,j,j')}\zeta_N^{l(-j+j'+1)}}{(\zeta_N^lz_{k,0}^{(N,j)})^{Nn-j+j'+1}}\\
        &=-\sum_{k=1}^{\infty}\frac{NC_{(N,j,j')}}{(z_{k,0}^{(N,j)})^{Nn-j+j'+1}}\\
        &=\begin{cases}
            \displaystyle-\frac{4C_{(4,0,j')}}{(\zeta_8\pi/\sqrt{2})^{4n+j'+1}}\sum_{k=1}^{\infty}\frac{1}{(2k-1)^{4n+j'+1}} & (\text{if }(N,j)=(4,0))\\
            \displaystyle-\frac{4C_{(4,2,j')}}{(\zeta_8\sqrt{2}\pi)^{4n-2+j'+1}}\sum_{k=1}^{\infty}\frac{1}{k^{4n-2+j'+1}} & (\text{if }(N,j)=(4,2))\\
            \displaystyle-\frac{6C_{(6,3,j')}}{(\zeta_{12}2\pi)^{6n-3+j'+1}}\sum_{k=1}^{\infty}\frac{1}{k^{6n-3+j'+1}} & (\text{if }(N,j)=(6,3))
        \end{cases}\\
        &=\begin{cases}
            \displaystyle\frac{4(-1)^n}{(\pi/\sqrt{2})^{4n+j'+1}}\sum_{k=1}^{\infty}\frac{1}{(2k-1)^{4n+j'+1}} & (\text{if }(N,j)=(4,0))\\
            \displaystyle\frac{4(-1)^{n+(j'-3)/2}}{(\sqrt{2}\pi)^{4n-2+j'+1}}\sum_{k=1}^{\infty}\frac{1}{k^{4n-2+j'+1}} & (\text{if }(N,j)=(4,2))\\
            \displaystyle\frac{6(-1)^{n+j'/2}}{(2\pi)^{6n-3+j'+1}}\sum_{k=1}^{\infty}\frac{1}{k^{6n-3+j'+1}} & (\text{if }(N,j)=(6,3))
        \end{cases}
    \end{align*}
    
    In the cases of $(N,j,j')=(4,0,3),\text{ }(4,0,1),\text{ }(4,2,3),\text{ and }(6,3,4)$, the condition $M=Nn-j+j'+1\ge\max\{2,-j+j'+1\}$ holds for any $n\ge0$, while in the other cases, the condition holds for any $n>0$. Replacing $n$ by $n-1$ in those four cases, we obtain \eqref{Equ_E40_and_lambda}, \eqref{Equ_E40_and_lambda_2}, \eqref{Equ_E42_and_zeta}, \eqref{Equ_E42_and_zeta_2}, \eqref{Equ_E63_and_zeta}, and \eqref{Equ_E63_and_zeta_2} for any positive integer $n$. The rest \eqref{Equ_E40_and_zeta} and \eqref{Equ_E40_and_zeta_2} is given by the relation $\zeta(s)=(2^s/(2^s-1))\times\lambda(s)$.

\end{proof}

\begin{corollary}
    There are the following equations between Bernoulli numbers $\{B_n\}$ and the congruential Euler numbers.
    \[B_{4n}=(-1)^n\frac{2n}{2^{2n}(2^{4n}-1)}\sum_{m=0}^{n-1}\binom{4n-1}{4m}\mathcal{E}_{4m}^{(4,0)}\quad(n>0)\]
    \[B_{4n-2}=(-1)^{n+1}\frac{4n-2}{2^{2n}(2^{4n-2}-1)}\sum_{m=0}^{n-1}\binom{4n-3}{4m}\mathcal{E}_{4m}^{(4,0)}\quad(n>0)\]
    \[B_{4n}=(-1)^n\frac{1}{2^{2n+1}(4n+1)}\sum_{m=0}^{n}\binom{4n+1}{4m}\mathcal{E}_{4m}^{(4,2)}\quad(n\ge0)\]
    \[B_{4n-2}=(-1)^{n+1}\frac{1}{2^{2n}(4n-1)}\sum_{m=0}^{n-1}\binom{4n-1}{4m}\mathcal{E}_{4m}^{(4,2)}\quad(n>0)\]
    \[B_{6n}=\frac{1}{3(6n+1)(6n+2)}\sum_{m=0}^{n}\binom{6n+2}{6m}\mathcal{E}_{6m}^{(6,3)}\quad(n\ge0)\]
    \[B_{6n-4}=\frac{1}{3(6n-2)(6n-3)}\sum_{m=0}^{n-1}\binom{6n-2}{6m}\mathcal{E}_{6m}^{(6,3)}\quad(n>0)\]
\end{corollary}
\begin{proof}
    Since $B_0=1$ and $\mathcal{E}_0^{(N,j)}=j!$, the third and fifth equations hold for $n=0$. For $n>0$, we use (\ref{Equ_E40_and_zeta}), (\ref{Equ_E40_and_zeta_2}), (\ref{Equ_E42_and_zeta}), (\ref{Equ_E42_and_zeta_2}), (\ref{Equ_E63_and_zeta}), (\ref{Equ_E63_and_zeta_2}), and a well known formula
    \[\zeta(2n)=\frac{(-1)^{n+1}(2\pi)^{2n}B_{2n}}{2(2n)!}.\]
\end{proof}

\appendix
\section{Background of the definition of congruential Euler numbers}

In 2019, in the paper \cite{Barman_and_Komatsu}, Barman and Komatsu introduced \textit{hypergeometric Lehmer-Euler numbers}\\$\{W_{N,n,r}^{(j)}\}$ by using the hypergeometric function in their generating function;
\begin{align*}
    \sum_{n=0}^{\infty}W_{N,n,r}^{(j)}\frac{z^n}{n!}&=\left({}_1F_r\left(1;\frac{rN+j+1}{r},\frac{rN+j+2}{r},\cdots,\frac{rN+j+r}{r};\left(\frac{t}{r}\right)^{r}\right)\right)^{-1}\\
    &=\left(\sum_{n=0}^{\infty}\frac{(rN+j)!}{(rN+rn+j)!}t^{rn}\right)^{-1}.
\end{align*}
This follows the original notation of \cite{Barman_and_Komatsu}. Note that the role of the parameter $N$ is different from that of the main text of this paper. Rather, $N$ of this paper corresponds to $r$ of Barman and Komatsu. In the following part of this appendix, $N$ refers to that of Barman and Komatsu.

They defined and considered hypergeometric Lehmer-Euler numbers for integers $N\ge0$ but only for $j=0,1$. In 2026, Kameyama, a former graduate student at Tohoku University, slightly changed the definition of hypergeometric Lehmer-Euler numbers of the case $N=0$ and extended them for the case $j>0$. He defined \textit{generalized $j$-Euler number} $\{\mathcal{E}_n^{(d,j)}\}$ by
\[\sum_{n=0}^{\infty}\mathcal{E}_{n}^{(d,j)}\frac{t^n}{n!}=\left(1+\sum_{l=1}^{\infty}\frac{t^{dl}}{(dl+j)!}\right)^{-1}\]
and researched their recurrence formula, explicit expressions, and determinant representations in his master's thesis \cite{Kameyama}(not available in open access). 

The original generating function of the hypergeometric Lehmer-Euler numbers (generalized j-Euler numbers) has the advantage of including a correction by $(rN+j)!$, which leads to the value of the constant term, i.e., $W_{N,0,r}^{(j)}$ being 1.
However, if $j\ge1$, then $\{W_{N,n,r}^{(j)}\}_{n\ge1}$ are not necessarily integers, regardless of whether the correction is applied. Rather, in the case $N=0$, defining it without correction allows for a simpler discussion and simpler representations of the results.

Therefore, in this paper, we slightly modify the definition of hypergeometric Lehmer-Euler numbers in the case $N=0$ by multiplying the original generating functions by $j!$.

This paper follows Kameyama's choice of character (that is $\mathcal{E}$) because, as shown in the paper, the relationship with the zeta values exhibits properties more similar to Euler numbers than to Lehmer numbers.

\section{Calculation data for the conjecture \ref{conjecture}}
In this appendix, we provide numerical evidence for Conjecture \ref{conjecture}.
Since the conjecture depends on three parameters $p,m,j$ and the amount of data becomes enormous, we present only representative examples for small values of $r$. All computations were carried out using SageMath.

\begin{conjecture}[Restatement of Conjecture \ref{conjecture}]
    Let $p$ be an prime number and $m$ be a divisor of $p-1$ or $m=2$. Let $q$ be a least common multiple of $2$ and $p-1$. For any positive integer $r$, there is a positive integer $n_0$ such that 
    \[\mathcal{E}_{mpn+qp^r}^{(mp,j)}\equiv\mathcal{E}_{mpn}^{(mp,j)}\mod{p^r}\quad\text{ for any }n\ge n_0.\]
\end{conjecture}

In the following, the ``period" means $P$ of $\mathcal{E}_{mpn}^{(mp,j)}\equiv\mathcal{E}_{mpn+P}^{(mp,j)}\mod{p^r}$. Thus, it is basically a multiple of $mp$.
\begin{longtable}{ccccc}
\caption{Numerical evidence for Conjecture \ref{conjecture}}\\
\toprule
Parameters & $p$ & $r$ & $n_0$ &Observed period \\
\midrule
$(mp,j)=(6,1)$ & $3$ & $1$ & $1$ & $6$ \\
$(mp,j)=(6,1)$ & $3$ & $2$ & $1$ & $18$ \\
$(mp,j)=(6,3)$ & $3$ & $1$ & $1$ & $6$ \\
$(mp,j)=(6,3)$ & $3$ & $2$ & $1$ & $18$ \\
$(mp,j)=(6,3)$ & $3$ & $3$ & $2$ & $54$ \\
$(mp,j)=(6,3)$ & $3$ & $4$ & $2$ & $162$ \\
$(mp,j)=(6,3)$ & $3$ & $5$ & $3$ & $486$ \\
$(mp,j)=(10,4)$ & $5$ & $1$ & $1$ & $20$ \\
$(mp,j)=(10,4)$ & $5$ & $2$ & $1$ & $100$ \\
$(mp,j)=(10,7)$ & $5$ & $3$ & $2$ & $500$ \\
$(mp,j)=(20,13)$ & $5$ & $1$ & $0$ & $20$ \\
$(mp,j)=(20,13)$ & $5$ & $2$ & $1$ & $100$ \\
$(mp,j)=(20,13)$ & $5$ & $3$ & $2$ & $500$ \\
$(mp,j)=(21,8)$ & $7$ & $1$ & $1$ & $42$ \\
$(mp,j)=(21,8)$ & $7$ & $2$ & $1$ & $294$ \\
$(mp,j)=(21,16)$ & $7$ & $1$ & $0$ & $21$ \\
$(mp,j)=(21,16)$ & $7$ & $2$ & $1$ & $294$ \\
$(mp,j)=(42,9)$ & $7$ & $1$ & $1$ & $21$ \\
$(mp,j)=(42,9)$ & $7$ & $2$ & $1$ & $294$ \\
$(mp,j)=(42,9)$ & $7$ & $2$ & $1$ & $2058$ \\
\bottomrule
\end{longtable}

The following is the concrete sequence data of the case of $(N,j)=(6,3)$. (As shown in section 4, $\mathcal{E}_{6n}^{(6,3)}$ is related to Bernoulli numbers.)

\begin{table}[htbp]
\centering
\begin{tabular}{ccc p{8cm}}
\toprule
$r$ & $n_0$ & periodic sequence \\
\midrule
$1$ & $1$ & $\overline{1}$ \\
$2$ & $1$ & $\overline{7,1,4}$ \\
$3$ & $2$ & $\overline{10,13,16,19,22,25,1,4,7}$ \\
$4$ & $2$ & $\overline{37,13,16,46,22,25,55,31,34,64,40,43,73,49,52,1,58,61,10,67,70,19,76,79,28,4,7}$ \\
\bottomrule
\end{tabular}
\caption{periodic sequence of $\mathcal{E}_{6n}^{(6,3)}$}
\end{table}

\end{document}